\documentclass[12pt]{amsart}
\usepackage[top=30truemm,bottom=30truemm,left=25truemm,right=25truemm]{geometry}
\usepackage{txfonts}
\usepackage{mathrsfs}
\usepackage{amsmath, amsthm, amssymb}
\usepackage{color}
\usepackage{bm}
\usepackage{amsfonts,amssymb}
\usepackage{dsfont}
\usepackage{amscd}
\usepackage{extarrows}
\usepackage{amsmath}
\usepackage{mathrsfs}
\usepackage{enumerate}
\usepackage{amscd}
\usepackage[all]{xy}
\usepackage{geometry}
\usepackage[colorlinks]{hyperref}
\usepackage[capitalize]{cleveref}
\newtheorem{theorem}{Theorem}[section]
\newtheorem{definition}{Definition}[section]
\newtheorem{lemma}{Lemma}[section]
\newtheorem{corollary}{Corollary}[section]
\newtheorem{proposition}{Proposition}[section]
\newtheorem{remark}{Remark}[section]

\newcommand{\PSH}[2]{\text{PSH}(#1, #2)}
\newcommand{\env}[2]{P_{#2}(#1)}

\newcommand{\be}{\begin{equation}}
	\newcommand{\ee}{\end{equation}}
\newcommand{\bea}{\begin{eqnarray}}
	\newcommand{\eea}{\end{eqnarray}}
\newcommand{\ben}{\begin{eqnarray*}}
	\newcommand{\een}{\end{eqnarray*}}
\newcommand{\bt}{\begin{split}}
	\newcommand{\et}{\end{split}}
\newcommand{\bet}{\begin{equation}}

\begin{document}
		\title[Complex Monge-Amp\`ere equations on compact Hermitian manifolds]
		{Degenerate complex Monge-Amp\`ere type equations on compact Hermitian manifolds and applications II}
		
		\author[H. Sun]{Haoyuan Sun}
		\address{Haoyuan Sun: School of Mathematical Sciences\\ Beijing Normal University\\ Beijing 100875\\ P. R. China}
		\email{202321130022@mail.bnu.edu.cn}
		\author[Z. Wang]{Zhiwei Wang}
		\address{Zhiwei Wang: Laboratory of Mathematics and Complex Systems (Ministry of Education)\\ School of Mathematical Sciences\\ Beijing Normal University\\ Beijing 100875\\ P. R. China}
		\email{zhiwei@bnu.edu.cn}
		
		\begin{abstract}
			Let $(X,\omega)$ be a compact Hermitian manifold of complex dimension $n$, equipped with a Hermitian metric $\omega$. Let $\beta$ be a possibly non-closed smooth $(1,1)$-form on $X$ such that $\int_X\beta^n>0$. Assume that there is a bounded $\beta$-plurisubharmonic function $\rho$ on $X$ and $\underline{\mathrm{Vol}}(\beta) > 0$. In this paper, we establish solutions to the degenerate complex Monge-Amp\`ere equations on $X$ within the Bott-Chern space of $\beta$ (as introduced by Boucksom-Guedj-Lu) and derive stability results for these solutions. As applications, we provide partial resolutions to the extended Tosatti-Weinkove conjecture and Demailly-P\u aun conjecture.
		\end{abstract}
		

		\maketitle
		
		\tableofcontents
		
		\section{Introduction}
		
		The complex Monge-Amp\`ere equation, which plays a fundamental role in pluripotential theory, complex analysis, and complex geometry, has led to many profound applications in K\"ahler geometry over the past decade. See, for example, \cite{Yau78,BT76,BT82,Ko98,Tsu88,Kol03,TZ06,EGZ09,BEGZ10,EGZ11,EGZ17} and the references therein.
		
		In recent years, there has been remarkable progress in Hermitian pluripotential theory, particularly in its applications to Hermitian geometry and degenerate complex Monge-Amp\`ere equations on Hermitian manifolds. See, for example, \cite{Che87,Han96,GL10,TW10,Kol05,KN15,Ngu16,LWZ24,BGL24} and the references therein.

		Let $(X, \omega)$ be a compact Hermitian manifold of complex dimension $n$, equipped with a Hermitian metric $\omega$. Let $\beta$ be a smooth possibly non-closed $(1,1)$ form. A function $u:X\rightarrow [-\infty,+\infty)$ is said to be quasi-plurisubharmonic if locally $u$ can be written as the sum of a smooth function and a plurisubharmonic function. A $\beta$-plurisubharmonic ($\beta$-psh for short) function $u$ is defined as a quasi-plurisubharmonic function satisfying $\beta+dd^cu\geq 0$ in the sense of currents. The set of all $\beta$-psh functions on $X$ is denoted by $\mbox{PSH}(X,\beta)$. We say $\beta$ is pseudoeffective, if $\mbox{PSH}(X,\beta)\neq \emptyset$. We say a pseudoeffective $(1,1)$-form is big (cf.\cite{BGL24}), if there is a $\beta$-psh function $\rho$ with analytic singularities such that $\beta+dd^c\rho$ dominates a Hermitian form on $X$.

		In \cite{LWZ24}, the authors solved the degenerate complex Monge-Amp\`ere equation
		\begin{equation*}
			(\beta + dd^c u)^n = f \, \omega^n
		\end{equation*}
		for $0 \leq f \in L^p(X)$ under the condition that the class $\beta$ is closed  and there is a bounded $\beta$-psh function (see also \cite{PSW24, LLZ24} for related results).
		
		In \cite{BGL24}, Boucksom-Guedj-Lu   made a remarkable extension of the foundational work by Boucksom-Eyssidieux-Guedj-Zeriahi \cite{BEGZ10} to the Hermitian setting.  Specifically, assume that  $\beta$ is a big (possibly non-closed) smooth $(1,1)$-form. They established that the degenerate complex Monge-Amp\`ere equation
		\begin{equation*}
			(\beta + dd^c \phi)^n = c f \omega^n
		\end{equation*}
		admits a solution for $0 \leq f \in L^p(X, \omega)$, where the constant $c > 0$ is uniquely determined by the data $(f, \beta, X, \omega, p, \|f\|_{p}:=\int_X|f|^p\omega^n)$.

		The main goal of this paper is to extend these results to possibly non-closed pseudoeffective $(1,1)$-forms $\beta$ that admit a bounded $\beta$-psh function, under the additional quantitative assumption  
		\begin{equation*}
			\underline{\mathrm{Vol}}(\beta):= \inf_{\substack{u \in \mathrm{PSH}(X,\beta) \cap L^\infty(X)}} \int_X (\beta + dd^c u)^n > 0.
		\end{equation*}
		Note that if $\beta$ is closed,  the condition $\underline{\mathrm{Vol}}(\beta)>0$ is equivalent to the condition that $\int_X\beta^n>0$, since $\underline{\mathrm{Vol}}(\beta)=\int_X\beta^n$.
		\begin{remark}
			Let $X$ be a compact complex surface, and let $\alpha$ be a semipositive, possibly non-closed $(1,1)$-form on $X$ satisfying $\int_X \alpha^2 > 0$ and $\partial\bar{\partial} \alpha = 0$ (for instance, $\alpha$ could be the Gauduchon metric on $X$). Let $\pi: \widetilde{X} \rightarrow X$ be a log resolution of $X$. Then, $\beta := \pi^* \alpha$ satisfies the desired assumptions. This observation can be utilized to construct higher-dimensional examples by considering products of compact complex surfaces with compact K\"ahler manifolds.
		\end{remark}
		
		Our first main result is as follows.
		\begin{theorem}[=\cref{MA equation 1}]\label{thm:main 1}
			Let $(X, \omega)$ be a compact Hermitian manifold of complex dimension $n$, equipped with a Hermitian metric $\omega$. Let $\beta$ be a possibly non-closed smooth $(1,1)$-form such that there exists a bounded $\beta$-psh function $\rho$ and $\underline{\mathrm{Vol}}(\beta) > 0$.
			Assume $0\leq f\in L^p(X,\omega^n)$ with $p>1$ and $||f||_p>0$. Then, for each $\lambda>0$, there exists a unique $\varphi_{\lambda}\in \mbox{PSH}(X,\beta)\cap L^{\infty}(X)$ such that
			$$
			(\beta+dd^c\varphi_{\lambda})^n=e^{\lambda\varphi_{\lambda}}f\omega^n.
			$$
			Furthermore, we have $|\varphi_{\lambda}|\leq C$ for some uniform constant C depending only on $\lambda,\beta,p,||f||_p$, $X,\omega$.
		\end{theorem}
		
		The key ingredient in the proof of \cref{thm:main 1} is the new $L^\infty$-a priori estimate established by Guedj-Lu in \cite{GL21}.
		When $\lambda=0$, we have the following theorem.
		\begin{theorem}[=\cref{main thm}]\label{thm:main 2}
			Let $(X, \omega)$ be a compact Hermitian manifold of complex dimension $n$, equipped with a Hermitian metric $\omega$. Let $\beta$ be a possibly non-closed smooth $(1,1)$-form such that there exists a bounded $\beta$-psh function $\rho$ and $\underline{\mathrm{Vol}}(\beta) > 0$.
			Assume that $0 \leq f \in L^p(X, \omega^n)$ with $p > 1$ and $\|f\|_p > 0$. Then there exists a constant $c > 0$ and a function $\varphi \in \mathrm{PSH}(X, \beta) \cap L^{\infty}(X)$ such that
			\[
			(\beta + dd^c\varphi)^n = cf \omega^n.
			\]
			Furthermore, the constant $c$ is uniquely determined by $f$, $\beta$, $X$ and the oscillation of $\varphi$, defined as $\mathrm{osc}\,\varphi := \sup_X \varphi - \inf_X \varphi$, satisfies $\mathrm{osc}\,\varphi \leq C$ for some constant $C$ depending only on $\beta$,  $p$, $\|f\|_{p}$, $X$, and $\omega$.
		\end{theorem}
		
		The proof of  \cref{thm:main 2} follows the strategy in \cite{GL23}. The uniqueness of the solution remains largely open. For the uniqueness when $f$ is uniformly bounded away from $0$ and $\beta$ is positive, see \cite{KN19}. 
		
		\begin{remark}
			When $\beta$ is closed, we have $\int_X \beta^n = \underline{\mathrm{Vol}}(\beta)$. In this case, the same results as those in Theorem \ref{thm:main 1} and Theorem \ref{thm:main 2} are obtained in \cite{LWZ24}. When $\beta$ is semipositive, the same results as those in Theorem \ref{thm:main 1} and Theorem \ref{thm:main 2} are obtained in \cite{GL23}.
		\end{remark}

		For applications, we establish the following stability result for the equation given in \cref{thm:main 1}.
		\begin{theorem}[=\cref{thm:stablity of MA 1}]\label{thm:main 3}
			
			Let $(X, \omega)$ be a compact Hermitian manifold of complex dimension $n$, equipped with a Hermitian metric $\omega$. Let $\beta$ be a possibly non-closed smooth $(1,1)$-form on $X$ such that there exists a bounded $\beta$-psh function $\rho$ and $\underline{\mathrm{Vol}}(\beta) > 0$. 
			Fix $f, g \in L^p(X, \omega^n)$ with $p > 1$ and $\|f\|_p, \|g\|_p>0$. Let $\varphi, \psi \in \mathrm{PSH}(X, \beta) \cap L^{\infty}(X)$ be such that
			\[
			(\beta + dd^c \varphi)^n = e^{\lambda \varphi} f \omega^n, \quad \text{and} \quad (\beta + dd^c \psi)^n = e^{\lambda \psi} g \omega^n,
			\]
			where $\lambda>0$ is a fixed constant. 
			Then we have the estimate
			\[
			\|\varphi - \psi\|_{\infty} \leq C \|f - g\|_p^{\frac{1}{n}},
			\]
			where $C > 0$ is a constant depending on $n$, $p$, $\lambda$, $\|f\|_p$, and $\|g\|_p$. 
			In particular, this implies the uniqueness of the solution to the complex  Monge-Amp\`ere equation in \cref{thm:main 1}.
		\end{theorem}	
		The above theorems can be used to establish a partial result toward an extended version of the Tosatti-Weinkove conjecture and the Demailly-P\u{a}un conjecture \cite[Conjectures 1.1, 1.2]{LWZ24}. Namely, we establish the following two theorems.
		
		\begin{theorem}[=\cref{thm: TW conj}]\label{thm: main 4}
			Let $(X, \omega)$ be a compact Hermitian manifold of complex dimension $n$, equipped with a Hermitian metric $\omega$. Let $\beta$ be a possibly non-closed smooth $(1,1)$-form on $X$ such that there exists a bounded $\beta$-psh function $\rho$ and $\underline{\mathrm{Vol}}(\beta) > 0$. Let $x_1,\dots,x_N$ be fixed points and $\tau_1,\dots,\tau_N$ be positive real numbers such that
			\begin{align*}
				\sum_{i=1}^{N}\tau_i^n < \underline{\mathrm{Vol}}(\beta). 
			\end{align*}
			Then there exists a $\beta$-psh function $\varphi$ with logarithmic poles at $x_1,\dots,x_N$:
			\begin{align*}
				\varphi(z) \leq c\tau_j\log|z| + O(1) 
			\end{align*}
			in a coordinate ball $(z_1,\dots,z_n)$ centered at $x_j$, where $c \geq 1$ is a uniform constant.
		\end{theorem}	
		
		\begin{definition}\label{def:upper lower volume}
			Let $(X, \omega)$ be a compact Hermitian manifold of complex dimension $n$, equipped with a Hermitian metric $\omega$. Let $\beta$ be a possibly non-closed smooth $(1,1)$-form on $X$ such that there exists a bounded $\beta$-psh function $\rho$. We then define:
			\begin{align*}
				\overline{\text{Vol}}_{n-1}(\beta) & := \sup_{u\in \text{PSH}(X, \beta)\cap L^{\infty}(X)}\int_X(\beta+dd^cu)^{n-1}\wedge\omega. 
			\end{align*}
			It is clear that the property $\overline{\text{Vol}}_{n-1}(\beta)<+\infty$ does not depend on the choice of the Hermitian metric $\omega$. 
		\end{definition}

		\begin{theorem}[=\cref{thm:DP conjecture}]\label{thm: main 5}
			Let $(X, \omega)$ be a compact Hermitian manifold of complex dimension $n$, equipped with a Hermitian metric $\omega$. Let $\beta$ be a possibly non-closed smooth $(1,1)$-form on $X$ such that there exists a bounded $\beta$-psh function $\rho$ and $\underline{\mathrm{Vol}}(\beta) > 0$  and $\overline{\mbox{Vol}}_{n-1}(\beta)<+\infty$.  Then $\beta$ is big, i.e., it contains a Hermitian current.
		\end{theorem}
		
		\begin{remark}
			In \cite[Theorem 3.20]{BGL24}, the authors proved that the function $\underline{\text{Vol}}$ vanishes outside the big cone under the assumption of the bounded mass property, i.e., the Hermitian metric $\omega$ satisfies
			\[
			\overline{\text{Vol}}(\omega) := \sup_{u \in \text{PSH}(X, \omega) \cap L^\infty(X)} \int_X (\omega + dd^c u)^n < +\infty.
			\]
			The bounded mass property directly implies that   $\overline{\text{Vol}}_{n-1}(\beta) < +\infty$.
			
			Indeed, it follows from \cite[Proposition 3.2]{GL22} that $\overline{\text{Vol}}(A\omega) < +\infty$ for all constants $A > 0$. Thus, we may assume that $\beta \leq \omega$.
			For any $u \in \text{PSH}(X, \beta) \cap L^\infty(X)$, $u \in \text{PSH}(X, \omega)$, and hence
			\[
			\int_X (\beta + dd^c u)^{n-1} \wedge \omega \leq \int_X (\omega + dd^c u)^{n-1} \wedge \omega \leq \int_X (\omega + dd^c u + \omega)^n \leq \overline{\text{Vol}}(2\omega) < +\infty.
			\]
			Consequently, our result extends \cite[Theorem 3.20]{BGL24} at least when the class $\{\beta\}$ admits a bounded $\beta$-psh potential.
		\end{remark}
		
		\begin{remark}
			When $\beta$ is closed, the same results as those in  \cref{thm: main 4} are obtained in \cite{LWZ24}. Additionally, when $\beta$ is closed and the Hermitian metric $\omega$ is pluriclosed, the same results as those in  \cref{thm: main 5} are obtained in \cite{LWZ24}.
		\end{remark}

		\subsection*{Acknowledgements}
		This research is supported by the National Key R\&D Program of China (Grant No. 2021YFA1002600). The second author is partially supported by grants from the National Natural Science Foundation of China (NSFC) (Grant No. 12071035) and by the Fundamental Research Funds for the Central Universities.

		\section{Preliminaries}
		Let $(X, \omega)$ be a compact Hermitian manifold of complex dimension $n$, equipped with a Hermitian metric $\omega$.
		
		\begin{definition}
			A function $u:X\rightarrow\mathbb{R}\cup\{-\infty\}$ is called quasi-plurisubharmonic if it can be locally written as a sum of a smooth function and a plurisubharmonic function.
			
			For a smooth $(1,1)$-form $\chi$ on $X$, a quasi-plurisubharmonic function $u$ is called $\chi$-psh, denoted $u\in\mbox{PSH}(X,\chi)$ if $\chi_u:=\chi+dd^cu\geq0$ in the weak sense of currents.
		\end{definition}
		
		Throughout this paper, we assume that $\beta$ is  a possibly non-closed smooth $(1,1)$-form such that there exists a bounded $\beta$-psh function $\rho$ and $\underline{\mathrm{Vol}}(\beta) > 0$.
		Following \cite{BGL24}, the Bott-chern space $BC^{p,q}(X)$ is defined as the cokernel of $dd^c:\Omega^{p-1,q-1}(X)\rightarrow\Omega^{p,q}(X)$, i.e.
		$$BC^{p,q}(X):=\frac{\Omega^{p,q}(X)}{dd^c\Omega^{p-1,q-1}(X)}.$$
		Note that for a bounded $\beta$-psh function $u$, the measure $(\beta+dd^cu)^n$ can be defined in the sense of Bedford-Taylor \cite{BT82}.
		
		We next establish the plurifine locality of bounded $\beta$-psh functions, originally due to \cite{BT87}.
		\begin{proposition}\label{maximum principle}
			Let $u,v$ be bounded $\beta$-psh functions on $X$, then we have
			$$
			\mathds{1}_{\{u>v\}}(\beta+dd^c\max(u,v))^n=\mathds{1}_{\{u>v\}}(\beta+dd^cu)^n.
			$$
			Furthermore,
			$$
			(\beta+dd^c\max(u,v))^n\geq\mathds{1}_{\{u>v\}}(\beta+dd^cu)^n+\mathds{1}_{\{u\leq v\}}(\beta+dd^cv)^n.
			$$
		\end{proposition}
		
		\begin{proof}
			
			We first check the first statement. Since it is a local property, we can verify it locally on an open subset $U$ of $X$. We can thus find a K\"ahler form $\eta \geq \beta$ on $U$. Then we have
			\begin{align*}
				&\mathds{1}_{\{u > v\}} (\beta + dd^c \max(u, v))^n \\
				= &\mathds{1}_{\{u > v\}} (\eta + dd^c \max(u, v) + \beta - \eta)^n \\
				= &\mathds{1}_{\{u > v\}} \sum_{k=1}^n \binom{n}{k} (\eta + dd^c \max(u, v))^k \wedge (\beta - \eta)^{n-k} \\
				= &\mathds{1}_{\{u > v\}} \sum_{k=1}^n \binom{n}{k} (\eta + dd^c u)^k \wedge (\beta - \eta)^{n-k} \\
				= &\mathds{1}_{\{u > v\}} (\beta + dd^c u)^n.
			\end{align*}
			The second equality follows from the standard Bedford-Taylor maximum principle.
			
			For the second statement, applying the first statement to the function $\max(u, v + \varepsilon)$, we get
			\[
			(\beta + dd^c \max(u, v + \varepsilon))^n \geq \mathds{1}_{\{u > v + \varepsilon\}} (\beta + dd^c u)^n + \mathds{1}_{\{u < v + \varepsilon\}} (\beta + dd^c v)^n.
			\]
			The conclusion follows by letting $\varepsilon \to 0$.
		\end{proof}
		\begin{corollary}\label{max principle2}
			For any bounded $\beta$-psh functions $u,v$ with $u\geq v$, we have
			$$
			\mathds{1}_{\{u=v\}}(\beta+dd^cu)^n\geq \mathds{1}_{\{u=v\}}(\beta+dd^cv)^n.
			$$
			
		\end{corollary}
		\begin{proof}
			This follows directly from the maximum principle Proposition \ref{maximum principle}.
		\end{proof}

		\subsection{Quasi-plurisubharmonic envelopes}
		If $h: X \to \mathbb{R} \cup \{\pm \infty\}$ is a measurable function, the $\beta$-plurisubharmonic (psh) envelope of $h$ is defined by
		$$
		\env{h}{\beta} = \left( \sup \left\{ \varphi \in \PSH{X}{\beta} : \varphi \leq h \quad \text{quasi-everywhere} \right\} \right)^*,
		$$
		with the convention that $\sup \emptyset = -\infty$. Here, quasi-everywhere means that $\varphi\leq h$ outside a pluripolar set of $X$. A subset $E\subset X$ is said to be pluripolar if $E$ is locally pluripolar, i.e., locally, $E$ is contained in the polar locus of a plurisubharmonic function.
		
		When $h = \min(u, v)$, we use the notation $\env{u, v}{\beta} := \env{\min(u, v)}{\beta}$ to denote the rooftop envelope of $u$ and $v$.
		
		\begin{lemma} \label{neg}
			Let $\{u_j\}_{j\in \mathbb{N}}$ be a sequence of $\beta$-psh function on $X$, which is  uniformly bounded above. Then the negligible set $\{ (\sup_{j} u_j)^* > \sup_{j} u_j\}$ is a pluripolar set. 
		\end{lemma}
		\begin{proof}
			
			The result is local; we may work on a small coordinate ball. On this ball, we can find a smooth function $\psi$ such that $\beta < dd^c \psi$ in the sense of smooth forms on $B$. It follows that $\psi + u_j \in \text{PSH}(B)$ for all $j$. Then we have
			\begin{equation*}
				\begin{split}
					&\left\{ \left(\sup_{j} u_j\right)^* > \sup_{j} u_j \right\} \\
					= &\left\{ \psi + \left(\sup_{j} u_j\right)^* > \psi + \sup_{j} u_j \right\} \\
					= &\left\{ \left(\sup_{j} (u_j + \psi)\right)^* > \sup_{j} (u_j + \psi) \right\}.
				\end{split}
			\end{equation*}
			The last set is negligible in the sense of \cite{BT82}, and thus, by \cite[Theorem 7.1]{BT82}, we deduce that the set above is pluripolar.
		\end{proof}
		
		\begin{corollary}[cf. {\cite[Page 7]{BGL24}}]\label{q.e.}
			For any measurable function $f:X\rightarrow[-\infty,+\infty]$ on $X$, one of the following three conditions holds:
			\begin{itemize}
				\item $P_{\beta}(f)\equiv-\infty$;
				\item $P_{\beta}(f)\equiv+\infty$;
				\item $P_{\beta}(f)$ is the largest $\beta$-psh function such that $P_{\beta}(f)\leq f$ quasi-everywhere.
			\end{itemize}
		\end{corollary}
		
		\begin{proof}
			If $P_{\beta}(f)$ is not identically $\infty$, it is of course $\beta$-psh by definition. By Lemma \ref{neg} we deduce that $P_{\beta}(f)\leq f$ quasi-everywhere. Moreover, the mean-value property shows easily the maximality of $P_{\beta}(f)$.
		\end{proof}
		
		\begin{remark}
			When $\beta$ is closed, we can prove that $P_{\beta}(f)<+\infty$ (i.e.,$P_{\beta}(f)$ is $\beta$-psh or identically $-\infty$) if and only if the set $\{f<+\infty\}$ is non-pluripolar, the argument is similar to that in \cite[Theorem 9.17]{GZ17}.
		\end{remark}

		\begin{theorem}\label{contact}
			Assume $f$ is quasi-continuous and bounded from below  on $X$, $P_{\beta}(f)\in \mbox{PSH}(X, \beta)$. Then $P_{\beta}(f)\leq f$ outside a pluripolar set and we have $$ \int_{\{ P_{\beta}(f)<f\}} (\beta+dd^c P_{\beta}(f))^n = 0.$$
			
		\end{theorem}
		\begin{proof}
			The first statement follows directly from Corollary \ref{q.e.}.
			
			Set $\beta_j := \beta + \varepsilon_j \omega_X$ for some decreasing sequence of positive numbers $\varepsilon_j \searrow 0$. Then $\beta_j$ is big, and it follows from \cite[Theorem 2.9]{BGL24} that for any $1 \leq j < k$,
			\[
			\int_{\{ P_{\beta_k}(f) < f \}} (\beta_j + dd^c P_{\beta_j}(f))^n = 0.
			\]
			Thanks to \cref{q.e.}, it is straightforward to check that $P_{\beta_j}(f)$ decreases to $P_{\beta}(f)$ as $j$ tends to $+\infty$. Since $f$ is assumed to be bounded from below and $\beta$ has a bounded potential, we have that $P_{\beta}(f)$ is bounded. The Bedford-Taylor convergence theorem tells us that the measure $(\beta_j + dd^c P_{\beta_j}(f))^n$ converges to $(\beta + dd^c P_{\beta}(f))^n$ weakly. Moreover, since the set $\{ P_{\beta_k}(f) < f \}$ is quasi-open, a standard argument gives that
			\[
			\int_{\{ P_{\beta_k}(f) < f \}} (\beta + dd^c P_{\beta}(f))^n = 0.
			\]
			Finally, we let $k \to \infty$ to conclude the proof.
		\end{proof}

		We will need the following two corollaries later.
		\begin{corollary}\label{contact2}
			Let $u,v$ be bounded $\beta$-psh functions and let $w:=P_{\beta}(u,v)=P_{\beta}(\min(u,v))$ be the rooftop envelope of $u,v$. Then
			$$
			(\beta+dd^cw)^n\leq\mathds{1}_{\{w=u<v\}}(\beta+dd^cu)^n+\mathds{1}_{\{w=v\}}(\beta+dd^cv)^n.
			$$
		\end{corollary}
		
		\begin{proof}
			Obviously, $\min(u, v)$ is quasi-continuous. Therefore, we can apply \cref{contact} to deduce that the measure $(\beta + dd^c w)^n$ is concentrated on the set $\{w = u\} \cup \{w = v\}$. The result then follows from \cref{max principle2}.
		\end{proof}
		
		\begin{corollary}\label{max principle 3}
			Fix $\lambda\geq0$ and $u,v\in \mbox{PSH}(X,\beta)\cap L^{\infty}(X)$, then the following holds:
			\begin{itemize}
				\item If $(\beta+dd^cu)^n\leq e^{\lambda u}f\omega^n$ and  $(\beta+dd^cv)^n\leq e^{\lambda v}g\omega^n$, then
				$$
				(\beta+dd^cP_{\beta}(u,v))^n\leq e^{\lambda P_{\beta}(u,v)}\max(f,g)\omega^n.
				$$
				\item If $(\beta+dd^cu)^n\geq e^{\lambda u}f\omega^n$ and  $(\beta+dd^cv)^n\geq e^{\lambda v}g\omega^n$, then
				$$
				(\beta+dd^c\max(u,v))^n\geq e^{\lambda \max(u,v)}\min(f,g)\omega^n.
				$$
			\end{itemize}
		\end{corollary}
		\begin{proof}
			The first and the second statements follow directly from \cref{contact2} and \cref{maximum principle}, respectively.
		\end{proof}

		\subsection{The domination principle}

		\begin{proposition}\label{domination for beta}
			Fix a constant $0\leq c<1$. If $u,v$ are bounded $\beta$-psh functions such that $\mathds{1}_{\{u<v\}}(\beta+dd^cu)^n\leq c\mathds{1}_{\{u<v\}}(\beta+dd^cv)^n$, then $u\geq v$.
		\end{proposition}
		
		\begin{proof}
			The proof follows from \cite[Proposition 2.8]{GL22}.
			
			Fix an arbitrary positive constant $\delta > 0$, we need to prove that $u \geq v - \delta$. Assume by contradiction that the set $E := \{u < v - \delta\}$ is non-empty and hence has positive Lebesgue measure since $u, v$ are quasi-psh. For $b > 1$, set
			\[
			u_b := P_{\beta}(bu - (b - 1)v).
			\]
			
			It follows from \cref{contact} that $(\beta + dd^c u_b)^n$ is concentrated on the contact set $D := \{u_b = bu - (b - 1)v\}$. Since $b^{-1}u_b + (1 - b^{-1})v \leq u$ with equality on $D$, by the maximum principle (\cref{maximum principle}) we get
			\begin{align*}
				&\mathds{1}_D(\beta + dd^c u)^n \geq \mathds{1}_D(\beta + dd^c(b^{-1}u_b + (1 - b^{-1})v))^n \\
				\geq{} & \mathds{1}_D b^{-n}(\beta + dd^c u_b)^n + \mathds{1}_D (1 - b^{-1})^n(\beta + dd^c v)^n \\
				\geq{} & \mathds{1}_D b^{-n}(\beta + dd^c u_b)^n + \mathds{1}_D c(\beta + dd^c v)^n,
			\end{align*}
			if we take $b$ sufficiently large. By our condition we have $\mathds{1}_{D \cap \{u < v\}}(\beta + dd^c u_b)^n = 0$.
			
			Since $u_b$ is bounded, by our assumption $\underline{\mbox{Vol}}(\beta) > 0$ we know the mass of $(\beta + dd^c u_b)^n$ on $D$ is positive and hence the set $D \cap \{u \geq v\}$ is non-empty. On this set,
			\[
			u_b = bu - (b - 1)v \geq u,
			\]
			thus the $\sup_X u_b$ is uniformly bounded from below since $u$ is assumed to be bounded. The sequence $u_b - \sup_X u_b$ converges in $L^1$ and almost everywhere to a function $u_{\infty} \in \mbox{PSH}(X, \beta)$ and hence $u_{\infty}$ must be identically $-\infty$ on $E$ with positive Lebesgue measure, this yields a contradiction that $u_{\infty} \notin \mbox{PSH}(X, \beta)$.
		\end{proof}

		\begin{corollary}
			\label{domination for beta 2}
			Let $u,v$ be bounded $\beta$-psh functions. If
			$$
			e^{-\lambda u}(\beta+dd^cu)^n\leq e^{-\lambda v}(\beta+dd^cv)^n,
			$$
			then $u\geq v$.
		\end{corollary}
		\begin{proof}
			For any fixed $\delta>0$, we have on the set $\{u<v-\delta\}$, we have
			$$
			(\beta+dd^cu)^n\leq e^{\lambda(u-v)}(\beta+dd^cv)^n \leq e^{-\lambda\delta}(\beta+dd^c(v-\delta))^n.
			$$
			It then follows from \cref{domination for beta} that $v-\delta\leq u$ and hence $v\leq u$ by letting $\delta\rightarrow0$.
		\end{proof}
		
		\begin{corollary}\label{domination for beta 3}
			Let $u,v$ be bounded $\beta$-psh functions. If
			$$
			(\beta+dd^cu)^n\leq c(\beta+dd^cv)^n
			$$
			for some constant $c>0$, then $c\geq 1$.
		\end{corollary}
		\begin{proof}
			Assume by contradiction that $c<1$. For any constant $C>0$, we have
			$$
			\mathds{1}_{\{u<v+C\}}(\beta+dd^cu)^n\leq c\mathds{1}_{\{u<v
				+C\}}(\beta+dd^c(v+C))^n.
			$$
			Using \cref{domination for beta} we deduce that $u\geq v+C$. Since $C$ was chosen arbitrarily and $u,v$ are bounded, this leads to a contradiction.
		\end{proof}
		
		\begin{remark}
			We remark that the domination principles above hold when $\beta$ satisfies
			$$
			\int_X(\beta+dd^cu)^n>0,
			$$
			for any $u\in\mbox{PSH}(X,\beta)\cap L^\infty(X)$. This condition is much weaker than the  condition
			$ \underline{Vol}(\beta)>0$.

		\end{remark}	
		\subsection{Mixed type inequalities}
		We recall the following result from \cite[Lemma 1.3]{GL23} (see also \cite[Lemma 1.9]{Ngu16}):
		
		\begin{lemma}\label{mixed type 1}
			Let $u_1,...,u_n$ be bounded $\alpha$-psh functions such that $(\alpha+dd^cu_i)^n= f_i\alpha^n$ for some $0\leq f_i\in L^1(X,\omega^n)$, then
			$$
			(\alpha+dd^cu_1)\wedge...\wedge(\alpha+dd^cu_n)\geq(\prod_{i=1}^nf_i)^{\frac{1}{n}}\alpha^n.
			$$
			Where $\alpha$ is a smooth semi-positive (1,1)-form on X such that $\int_X\alpha^n>0$.
		\end{lemma}
		
		We will need the following more general version of mixed type inequalities:
		
		\begin{lemma}\label{mixed type 2}
			Let $\beta_1,...,\beta_n$ be smooth $(1,1)$-forms admitting bounded potentials $\rho_1\in \mbox{PSH}(X,\beta_1)\cap L^\infty(X),...,\rho_n\in \mbox{PSH}(X,\beta_n)\cap L^\infty(X)$. For $1\leq j\leq n$, let $u_j$ be bounded $\beta_j$-psh functions such that $(\beta_j+dd^cu_j)^n=f_j\omega^n$ for some $0\leq f_j\in L^1(X,\omega^n)$, then
			$$
			(\beta_1+dd^cu_1)\wedge...\wedge(\beta_n+dd^cu_n)\geq(\prod_{j=1}^nf_j)^{\frac{1}{n}}\omega^n.
			$$
			Where $\omega$ is a fixed Hermitian metric.
		\end{lemma}
		We need three preliminary lemmas before going to prove \cref{mixed type 2}.

		\begin{lemma}\label{local dirichlet for chi}
			Let $0 \leq f \in L^p(\Omega)$ for some $p > 1$ and let $\chi$ be an arbitrary smooth $(1,1)$-form defined near $\overline{\Omega}$, where $\Omega$ is a bounded strictly pseudoconvex domain in $\mathbb{C}^n$. Assume $\varphi \in C^0(\partial\Omega)$. Then there exists a unique continuous function $u \in \mathrm{PSH}(\Omega,\chi) \cap C^0(\overline{\Omega})$ that solves
			\[
			(\chi + dd^c u)^n = f\omega^n,\quad u = \varphi \text{ on } \partial\Omega.
			\]
		\end{lemma}
		
		\begin{proof}
			Let $\eta$ be the strictly plurisubharmonic defining function of $\Omega$ (which is defined near $\overline{\Omega}$). Up to rescaling, we may assume that $\chi_1 := \chi + dd^c\eta > 0$ is a Hermitian form on $\Omega$. Using \cite[Theorem 4.2]{KN15}, we see that we can solve the following Dirichlet problem:
			\[
			\begin{cases}
				(\chi_1 + dd^c u)^n = f\omega^n, & \text{on } \Omega, \\
				u = \varphi - \eta, & \text{on } \partial\Omega, \\
				u\in C(\overline{\Omega})\cap \mathrm{PSH}(\Omega,\chi_1).
			\end{cases}
			\]
			It follows immediately that $u + \eta$ is the desired solution.
		\end{proof}

		\begin{lemma}\label{local stability for chi}
			Let $\Omega$ be a bounded strictly pseudoconvex domain in $\mathbb{C}^n$, let $\chi$ be an arbitrary smooth $(1,1)$-form defined near $\overline{\Omega}$ and let $f_1,f_2$ be non-negative functions in $L^p(\Omega)$, $p>1$. Let $\phi_1,\phi_2\in C^0(\partial\Omega)$. Suppose $u_1,u_2$ solve the following Dirichlet problems:
			\[
			\begin{cases}
				(\chi + dd^c u_j)^n = f_j \omega^n, & \text{on } \Omega, \\
				u_j = \phi_j, & \text{on } \partial\Omega, \\
				u_j \in C(\overline{\Omega}) \cap \mathrm{PSH}(\Omega,\chi)
			\end{cases}
			\]
			for $j = 1, 2$. Then,
			\[
			\|u_1 - u_2\|_{L^{\infty}(\Omega)} \leq \sup_{\partial\Omega} |\phi_1 - \phi_2| + C \|f_1 - f_2\|_{L^p(\Omega)}^\frac{1}{m},
			\]
			where $C$ depends only on $p$ and $\Omega$.
		\end{lemma}
		
		\begin{proof}
			A similar argument as in \cref{local dirichlet for chi} transforms the problem to the case of \cite[Theorem 4.1]{KN15}.
		\end{proof}
		\begin{lemma}\label{local comparison principle}
			Let $\Omega$ be a bounded strictly pseudoconvex domain in $\mathbb{C}^n$, let $\chi$ be an arbitrary smooth $(1,1)$-form defined near $\overline{\Omega}$ and let $u,v\in \mbox{PSH}(X,\chi)\cap L^{\infty}(\Omega)$ be such that $\underset{\zeta\rightarrow z\in\partial\omega}{\liminf}(u-v)(\zeta)\geq0$. Suppose that $(\chi+dd^cu)^n\leq(\chi+dd^cv)^n$, then $v\leq u$ on $\Omega$.
		\end{lemma}	
		\begin{proof}
			The same argument as in \cref{local dirichlet for chi} transforms the problem to the case of \cite[Corollary 3.4]{KN15}.
		\end{proof}
		\begin{proof}[Proof of \cref{mixed type 2}]
			The arguments is adapted from \cite[Lemma 6.3]{Kol05}. The problem is local, so we can work on a small coordinate ball $B\subset X$.
			We first assume all the $u_k$ are continuous near $\overline{B}$, choose sequences of smooth positive functions $f_k^j\rightarrow f_k$ in $L^p(\overline{B},\omega^n)$ and smooth functions $\phi_k^j$ converge uniformly to $u_k$ on $\partial B$ for each $k$. Since we are working on a ball $B$, the subsolution condition in \cite[Theorem 1.1]{GL10} is automatically satisfied, thus we can find sequences $v_k^j\in \mbox{PSH}(B,\beta_k)\cap C^{\infty}(X)$ solving the following system of Dirichlet problems:
			$$
			\begin{cases}
				(\beta_k+dd^cv_k^j)^n=f_k^j\omega^n,  & on\quad B \\
				v_k^j=\phi_k^j,    & on\quad\partial B
			\end{cases}.
			$$
			By the stability estimate \cref{local stability for chi},
			$$
			||u_k-v_k^j||_{\infty}\leq \underset{\partial B}{\sup}\,|u_k-v_k^j|+C(p,n,B)||f_k-f_k^j||_{L^p}^{\frac{1}{n}}.
			$$
			and hence $v_k^j\rightarrow u_k$ uniformly as $j\rightarrow\infty$ for all $k$. Due to the concavity of the mapping $A \mapsto \log\left(\det^{1/n}A\right)$ defined on the set of positive definite Hermitian matrices, and by using a limit argument, we have
			\[
			(\beta_1 + dd^c v_1^j) \wedge \cdots \wedge (\beta_n + dd^c v_n^j) \geq \left(\prod_{k = 1}^n f_k^j\right)^{\frac{1}{n}} \omega^n.
			\]
			Letting $j \to \infty$, we obtain the desired inequality.
			
			Next, we address the case when all the $f_k$ belong to $L^p(X,\omega^n)$ for some $p>1$. Since $u_k$ is upper semi-continuous, we can find sequences of continuous functions $\psi_k^j$ decreasing to $u_k$ on $\overline{B}$ for each $k$. Using \cref{local dirichlet for chi}, we can solve the following system of Dirichlet problems:
			$$
			\begin{cases}
				(\beta_k+dd^cu_k^j)^n=f_k\omega^n,  & on\quad B; \\
				u_k^j=\psi_k^j,    & on\quad\partial B;\\
				u_k^j\in C(\overline{B})\cap \mbox{PSH}(B,\beta_k).
			\end{cases}
			$$
			Since $\psi_k^j$ decreases to $u_k$,  it follows from the comparison principle \cref{local comparison principle} that $u_k^j$ is decreasing with respect to $j$ for each $k$ and $u_k^j\geq u_k$. Then it is clear that $u_k^j$ decreases to $u_k$. Thus the claim follows from the decreasing convergence theorem of Bedford-Taylor \cite{BT82}.
			
			Finally, for the general case, we can approximate $f_k$ by a sequence of increasing  $L^2$ functions and proceed similarly as in the previous steps, as described in \cite[Lemma 6.3]{Kol05}. The proof is thus complete.
		\end{proof}

		\section{$L^{\infty}$ a priori estimates}
		
		In this section, we establish the $L^{\infty}$ a priori estimates for Monge-Amp\`ere equations in our setting. We will follow a method  in \cite{GL21} and \cite{GL23} (see also \cite{DDL25} for the relative case and \cite{PSW24} when $\beta$ is closed).
		
		\begin{lemma}\label{lem: concave weight}Let $(X,\omega)$ be a compact Hermitian manifold of complex dimension $n$, equipped with a Hermitian metric $\omega$. 
			Let $\beta$ be  a possibly non-closed smooth $(1,1)$-form such that there exists a bounded $\beta$-psh function $\rho$ and $\underline{\text{Vol}}(\beta) > 0$.
			Fix a concave increasing function $\chi:\mathbb{R}^-\rightarrow\mathbb{R}^-$ such that $\chi'(0)\geq 1$. Let $\varphi\leq\rho$ be a bounded $\beta$-psh function and $\psi=\chi\circ(\varphi-\rho)+\rho$. Then
			$$
			(\beta+dd^cP_{\beta}(\psi))^n\leq\mathds{1}_{\{P_{\beta}(\psi)=\psi\}}(\chi'\circ(\varphi-\rho))^n(\beta+dd^c\varphi)^n.
			$$
		\end{lemma}
		
		\begin{proof}
			The proof follows the argument in \cite[Lemma 1.6]{GL21}. Up to approximation we can assume without loss of generality that $\chi$ is smooth. 
			Set $\sigma:=\chi^{-1}:\mathbb{R}^-\rightarrow\mathbb{R}^-$ be the inverse function of $\chi$, it is easy to see that $\sigma$ is a convex increasing function such that $\sigma'\leq 1$. Let $\eta:=P_{\beta}(\psi)-\rho$ and $v:=\sigma\circ\eta+\rho$.
			\begin{align*}
				\beta+dd^cv=&\beta+dd^c{\rho}+\sigma''\circ\eta d\eta\wedge d^c\eta+\sigma'\circ\eta dd^c\eta\\
				\geq&(1-\sigma'\circ\eta)(\beta+dd^c{\rho})+\sigma'\circ\eta(\beta+dd^cP_{\beta}(\psi))\\
				\geq&\sigma'\circ\eta(\beta+dd^cP_{\beta}(\psi))\geq0.
			\end{align*}
			It follows that $v$ is $\beta$-psh and $(\beta+dd^cP_{\beta}(\psi))^n\leq1_{\{P_{\beta}(\psi)=\psi\}}(\sigma'\circ\eta)^{-n}(\beta+dd^cv)^n$. On the contact set $\{P_{\beta}(\psi)=\psi\}$ we have 
			$$
			\sigma'\circ\eta=\sigma'\circ(\psi-\rho)=\frac{1}{(\chi'\circ(\psi-\rho))}.
			$$
			We also have $v\leq\sigma\circ(\psi-\rho)+\rho=\varphi$ quasi everywhere, hence everywhere on $X$ since both functions are quasi-psh. By \cref{max principle2} we obtain that $1_{\{P_{\beta}(\psi)=\psi\}}(\beta+dd^cv)^n\leq1_{\{P_{\beta}(\psi)=\psi\}}(\beta+dd^c\varphi)^n$ and conclude the proof.
		\end{proof}
		
		\begin{theorem}\label{thm: a prior estimate}
			Let $\beta,\rho$ be as above. Let $\mu$ be a probability measure on X such that $\mbox{PSH}(X,\beta)\subset L^m(\mu)$ for some $m>n$. Then any solution $\varphi\in  \mbox{PSH}(X,\beta)\cap L^{\infty}(X)$ to $(\beta+dd^c\varphi)^n\leq c\mu$, where $c>0$ is an appropriate constant, satisfies
			$$
			osc_X(\varphi)\leq T
			$$
			for some uniform constant T which only depends on $X$, the bound of $\rho$, the upper bound of $\frac{c}{\underline{vol}(\beta)}$ and
			$$
			A_m(\mu):=\sup\left\{\left(\int_X(-\psi)^md\mu\right)^{\frac{1}{m}}: \psi\in \mbox{PSH}(X,\beta)\, \mbox{with} \,\underset{X}{\sup}\,\psi=0\right\}.
			$$ 
		\end{theorem}
		
		\begin{remark}
			As noted in \cite{GL21}, this theorem applies to measures $\mu=f\omega^n$ for $0\leq f\in L^p(X)$ with $p>1$. It also applies to more general densities, such as the Orlicz weight. We refer the reader to \cite{GL21} and \cite{GL23} for more details.
		\end{remark}
		
		\begin{proof}[Proof of \cref{thm: a prior estimate}]
			The proof follows from \cite[Theorem 2.2]{GL23}, with slight modifications in the relative case. 
			We define $\rho$ as the $\beta$-psh envelope of all non-positive $\beta$-psh functions. Without loss of generality, we may assume that
			\[
			\sup_X \varphi = 0.
			\]
			It is clear that $\varphi \leq \rho$. Now, set
			\[
			T_{\text{max}} := \sup \{ t > 0 : \mu(\{\varphi < \rho - t\}) > 0 \}.
			\]
			By the definition of $T_{\text{max}}$, we have
			\[
			\mu(\{\varphi < \rho - T_{\text{max}}\}) = 0.
			\]
			This implies that $\varphi \geq \rho - T_{\text{max}}$ $\mu$-almost everywhere.
			Since we have $c\mu\geq(\beta+dd^c\varphi)^n$, the domination principle \cref{domination for beta}  (take $c=0$, $u=\varphi$, $v=\rho-T_{max}$) implies that $\varphi\geq\rho-T_{max}$. Thus we only need to establish the bound for $T_{max}$.
			
			Let $\chi:\mathbb{R}^-\rightarrow\mathbb{R}^-$ be a concave increasing weight such that $\chi(0)=0$ and $\chi'(0)=1$. Set $\psi:=\chi\circ(\varphi-\rho)+\rho$ and $u:=P_{\beta}(\psi)$. By \cref{lem: concave weight} we have
			$$
			\frac{1}{c}(\beta+dd^cu)^n\leq1_{\{u=\psi\}}(\chi'\circ(\varphi-\rho))^n\mu.
			$$
			\textbf{Step 1.} We begin by controlling the energy of the function $u$. Let $\varepsilon > 0$ be a small positive constant, and define $m := n + 3\varepsilon$. Given that $\chi(0) = 0$ and $\chi$ is concave, we have the inequality $|\chi(t)| \leq |t|\chi'(t)$ for all $t$. We can now estimate the energy as follows:
			
			\begin{align*}
				\int_X (\rho - u)^{\varepsilon} \frac{(\beta + dd^c u)^n}{c} \leq{} & \int_X (-\chi \circ (\varphi - \rho))^{\varepsilon} (\chi' \circ (\varphi - \rho))^n d\mu \\
				\leq{} & \int_X (\rho - \varphi)^{\varepsilon} (\chi' \circ (\varphi - \rho))^{n + \varepsilon} d\mu \\
				\leq{} & \left(\int_X (\rho - \varphi)^{n + 2\varepsilon} d\mu\right)^{\frac{\varepsilon}{n + 2\varepsilon}} \left(\int_X (\chi' \circ (\varphi - \rho))^{n + 2\varepsilon} d\mu\right)^{\frac{n + \varepsilon}{n + 2\varepsilon}} \\
				\leq{} & \left(\int_X (-\varphi)^{n + 2\varepsilon} d\mu\right)^{\frac{\varepsilon}{n + 2\varepsilon}} \left(\int_X (\chi' \circ (\varphi - \rho))^{n + 2\varepsilon} d\mu\right)^{\frac{n + \varepsilon}{n + 2\varepsilon}} \\
				\leq{} & A_m(\mu)^{\varepsilon} \left(\int_X (\chi' \circ (\varphi - \rho))^{n + 2\varepsilon} d\mu\right)^{\frac{n + \varepsilon}{n + 2\varepsilon}}.
			\end{align*}
			Here, the first inequality follows from the fact that $u = \chi \circ (\varphi - \rho) + \rho$ on the support of $(\beta + dd^c u)^n$. The second inequality is derived by applying the inequality $|\chi(t)| \leq |t|\chi'(t)$. The third inequality is a consequence of H\"older's inequality. The last inequality is the definition of $A_m(\mu)$.
			
			\textbf{Step 2.} We select an appropriate weight function $\chi$ such that the integral 
			\[
			\int_X \left(\chi' \circ (\varphi - \rho)\right)^{n + 2\varepsilon} \, d\mu \leq 2
			\]
			is under control.
			
			To begin with, we fix a value $T \in [0, T_{\max})$. Suppose $g: \mathbb{R}^+ \rightarrow \mathbb{R}^+$ is an increasing, absolutely continuous function on the interval $[0, T]$ with $g(0) = 1$. Then, we have the following identity:
			
			\[
			\int_X g \circ (\rho - \varphi) \, d\mu = \mu(X) + \int_0^{T_{\max}} g'(t) \mu(\varphi < \rho - t) \, dt.
			\]
			
			Next, we define the function $f(t)$ as follows:
			
			\[
			f(t) := 
			\begin{cases} 
				\dfrac{1}{(1 + t)^2 \mu(\varphi < \rho - t)}, & \text{if } t \in [0, T]; \\
				\dfrac{1}{(1 + t)^2}, & \text{if } t > T. 
			\end{cases}
			\]
			
			It can be verified that $f(t)$ is an $L^1$ function. We then define $g(x)$ as:
			
			\[
			g(x) := \int_0^x f(t) \, dt + 1.
			\]
			
			It is straightforward to check that $g$ is absolutely continuous on $[0, T]$ and that $g'(t) = f(t)$ almost everywhere (specifically, at every Lebesgue point).
			
			Now, we redefine $-\chi(-x)$ as:
			
			\[
			-\chi(-x) := \int_0^x g(t)^{\frac{1}{n + 2\varepsilon}} \, dt.
			\]
			
			From this definition, we obtain $\chi(0) = 0$, $\chi'(0) = 1$, and $g(t) = \left(\chi'(-t)\right)^{n + 2\varepsilon}$. Note that this equality holds everywhere because $g$ is continuous and $g \geq 1$. This construction ensures that $\chi$ is a concave, increasing function with $\chi' \geq 1$. Consequently, we have:
			
			\[
			\int_X \left(\chi' \circ (\varphi - \rho)\right)^{n + 2\varepsilon} \, d\mu \leq \mu(X) + \int_0^{+\infty} \frac{1}{(1 + t)^2} \, dt \leq 2.
			\]

			\textbf{Step 3.}
			Next, we estimate the level set $\mu(\varphi < \rho - t)$. By utilizing Step 2, we can establish a uniform lower bound for $\sup_{X} u$ as follows:
			
			\begin{align*}
				\left(-\sup_{X}(u - \rho)\right)^{\varepsilon} \frac{\underline{\mathrm{Vol}}(\beta)}{c} &\leq \left(-\sup_{X}(u - \rho)\right)^{\varepsilon} \int_X \frac{(\beta + dd^c u)^n}{c} \\
				&\leq \int_X (\rho - u)^{\varepsilon} \frac{(\beta + dd^c u)^n}{c} \\
				&\leq 2 A_m(\mu)^{\varepsilon}.
			\end{align*}
			
			This implies that $0 \geq \sup_{X} (u - \rho) \geq -2^{\frac{1}{\varepsilon}} \left(\frac{c}{\underline{\mathrm{Vol}}(\beta)}\right)^{\frac{1}{\varepsilon}} A_m(\mu)$. From this, we infer that
			
			\[
			\|\chi \circ (\varphi - \rho)\|_{L^m(\mu)} \leq \|u - \rho\|_{L^m(\mu)} \leq \left(1 + 2^{\frac{1}{\varepsilon}} \left(\frac{c}{\underline{\mathrm{Vol}}(\beta)}\right)^{\frac{1}{\varepsilon}}\right) A_m(\mu) := B.
			\]
			
			Indeed, there exists a constant $c$ (which may depend on other parameters) such that $c \leq 2^{\frac{1}{\varepsilon}} \left(\frac{c}{\underline{\mathrm{Vol}}(\beta)}\right)^{\frac{1}{\varepsilon}} A_m(\mu)$ and $\sup_{X}(u - \rho) + c = 0$. By definition, we have $\|u - \rho + c\|_{L^m(\mu)} \leq A_m(\mu)$, and thus, by the triangle inequality, $\|u - \rho\|_{L^m(\mu)} \leq c + A_m(\mu) \leq \left(1 + 2^{\frac{1}{\varepsilon}} \left(\frac{c}{\underline{\mathrm{Vol}}(\beta)}\right)^{\frac{1}{\varepsilon}}\right) A_m(\mu)$.
			
			Given that $u - \rho \leq \chi \circ (\varphi - \rho) \leq 0$ quasi everywhere (note that $\mu$ does not charge pluripolar sets under our assumptions), we deduce that $\|\chi \circ (\varphi - \rho)\|_{L^m(\mu)} \leq \|u - \rho\|_{L^m(\mu)} \leq B$. Consequently, we have
			
			\begin{align*}
				\mu(\varphi < \rho - t) &= \int_{\{\varphi < \rho - t\}} d\mu \\
				&\leq \int_X \frac{|\chi \circ (\varphi - \rho)|^m}{|\chi(-t)|^m} d\mu \\
				&= \frac{1}{|\chi(-t)|^m} \|\chi \circ (\varphi - \rho)\|_{L^m(\mu)}^m \\
				&\leq \frac{B^m}{|\chi(-t)|^m}.
			\end{align*}

			\textbf{Step 4.} Conclusion.  Let us define $h(t) := -\chi(-t)$. We observe that $h(0) = 0$ and $h'(t) = [g(t)]^{\frac{1}{n + 2\varepsilon}}$ is positive and increasing, implying that $h$ is convex. Additionally, since $g(t) \geq g(0) \geq 1$, we have $h'(t) \geq 1$ and $h(1) = \int_0^1 h'(t) \, dt \geq 1$.
			
			For almost all $t \in [0, T]$, we have the inequality:
			
			\[
			\frac{1}{(1 + t)^2 g'(t)} = \mu(\varphi < \rho - t) \leq \frac{B^m}{h^m(t)}.
			\]
			
			This leads to:
			
			\[
			h^m(t) \leq B^m (1 + t)^2 g'(t) = (n + 2\varepsilon) B^m (1 + t)^2 h''(t) (h')^{n + 2\varepsilon - 1}(t),
			\]
			which holds almost everywhere.
			
			Next, we multiply both sides by $h'$ and integrate from $0$ to $t$:
			
			\begin{align*}
				\frac{h^{m + 1}(t)}{m + 1} &\leq (n + 2\varepsilon) B^m \int_0^t h''(s) (h')^{n + 2\varepsilon}(s) \, ds \\
				&\leq \frac{(n + 2\varepsilon) B^m (1 + t)^2}{n + 2\varepsilon + 1} \left( (h')^{n + 2\varepsilon + 1}(t) - 1 \right) \\
				&\leq B^m (1 + t)^2 (h')^{n + 2\varepsilon + 1}(t).
			\end{align*}
			The first and second inequalities hold because $h$ and $h'$ are absolutely continuous, allowing us to take derivatives and integrate.
			
			Recall that $\alpha := m + 1 = n + 3\varepsilon + 1 > \gamma := n + 2\varepsilon + 1 > 2$. The previous inequality can be rewritten as:
			
			\[
			(1 + t)^{-\frac{2}{\gamma}} \leq C h'(t) h(t)^{-\frac{\alpha}{\gamma}},
			\]
			for some uniform constant $C$ depending only on $n$, $m$, and $B$. Integrating both sides from $1$ to $T$, we obtain:
			
			\[
			\frac{(1 + T)^{1 - \frac{2}{\gamma}}}{1 - \frac{2}{\gamma}} - \frac{2^{1 - \frac{2}{\gamma}}}{1 - \frac{2}{\gamma}} \leq C \int_1^T \frac{(h^{1 - \frac{\alpha}{\gamma}})'}{1 - \frac{\alpha}{\gamma}} \leq C_1 h(1)^{1 - \frac{\alpha}{\gamma}} \leq C_1,
			\]
			since $1 - \frac{\alpha}{\gamma} < 0$ and $h(1) \geq 1$. From this, we deduce that $T$ can be bounded by a constant $C_2$ that depends only on $A_m(\mu)$.
			
			Finally, by letting $T \rightarrow T_{\max}$, we conclude the proof.
		\end{proof}
		
		\section{Construction of the subsolution} 
		In this section, we prove the following lemma establishing a subsolution to the complex Monge-Ampère equation, adapting the approach from \cite[Lemma 3.3]{GL23}.
		\begin{lemma}\label{subsolution} Let $(X,\omega)$ be a compact Hermitian manifold of complex dimension $n$, equipped with a Hermitian metric $\omega$. 
			Let $\beta$ be  a possibly non-closed smooth $(1,1)$-form such that there exists a bounded $\beta$-psh function $\rho$ and $\underline{\text{Vol}}(\beta) > 0$.
			Fix a constant $p>1$. Then there exists a constant $c>0$ depending only on $\beta$ and $p$ such that for any $0\leq f\in L^p(\omega^n)$ with $||f||_{p}\leq1$, we can find $u\in \mbox{PSH}(X,\beta)\cap L^{\infty}(X)$ such that $-T\leq u\leq0$ and
			$$
			(\beta+dd^cu)^n\geq cf\omega^n.
			$$
			Where $T$ is the uniform a priori bound in \cref{thm: a prior estimate}.
		\end{lemma}
		\begin{proof}
			
			Let $\{\beta_j\}$ be defined as $\beta_j := \beta + \varepsilon_j \omega$ for some sequence $\{\varepsilon_j\}$ decreasing to $0$. Without loss of generality, we may assume that $\beta < \beta_j \leq \omega$ by rescaling $\omega$ if necessary. Since $\beta_j$ is a Hermitian class, we can apply \cite[Theorem 0.1]{KN15} (see also \cite[Theorem 4.6]{BGL24}) to solve the equation
			\[
			(\beta_j + dd^c u_j)^n = c_j (1 + f) \omega^n, \quad \text{with} \quad \sup_X u_j = 0,
			\]
			where $u_j \in \text{PSH}(X, \beta_j) \cap L^{\infty}(X)$ and $c_j > 0$.
			
			Let $a$ be a fixed positive constant such that $a < 1$. By our assumption, $(\omega + dd^c u_j)^n \geq (\beta_j + dd^c u_j)^n \geq c_j \omega^n$. Hence, we can use the mixed-type inequality (denoted as \cref{mixed type 1} in the context of the paper) to deduce that
			\[
			(\omega + dd^c u_j) \wedge \omega^{n-1} \geq a^{\frac{n-1}{n}} c_j^{\frac{1}{n}} \omega^n.
			\]
			Integrating both sides, we obtain
			\begin{align*}
				a^{\frac{n-1}{n}} c_j^{\frac{1}{n}} \int_X \omega^n &\leq \int_X (\omega + dd^c u_j) \wedge \omega^{n-1} \\
				&= \int_X \omega^n + \left| \int_X u_j dd^c \omega^{n-1} \right| \\
				&\leq \int_X \omega^n + B \int_X (-u_j) \omega^n \\
				&\leq C,
			\end{align*}
			where in the last inequality, we have chosen a constant $B$ such that $dd^c \omega^{n-1} \leq B \omega^n$ and used the compactness of psh functions (see, e.g., \cite[Proposition 1.1]{Ngu16}). From the above, we conclude that the sequence $\{c_j\}$ is uniformly bounded from above.
			
			We next claim that $\{u_j\}$ is uniformly bounded. Since we can choose a common $\rho$ for all $\beta_j$ and all $A_m(\mu, \beta_j)$ are bounded by $A_m(\mu, \beta_1)$, by \cref{thm: a prior estimate}, it suffices to show that $\underline{\text{Vol}}(\beta_j)$ is uniformly bounded below by a positive constant. However, we have
			\[
			\underline{\text{Vol}}(\beta_j) = \underline{\text{Vol}}(\beta_j + dd^c \rho) \geq \underline{\text{Vol}}(\beta + dd^c \rho) = \underline{\text{Vol}}(\beta) > 0,
			\]
			where the first and third equalities follow from the definition of the lower volume of a positive current as in \cite{BGL24}, and the second inequality follows from \cite[Proposition 3.7]{BGL24}. The claim follows.
			
			Up to extracting a subsequence, we may assume that $u_j \rightarrow u \in \text{PSH}(X, \beta)$, and \cref{thm: a prior estimate} shows that $-T \leq u \leq 0$ for some uniform constant $T > 0$. Set
			\[
			\Phi_j := \left( \sup_{k \geq j} u_k \right)^*.
			\]
			By Hartogs' lemma, we have $\Phi_j$ decreases to $u$. Recall that we have
			\[
			(\beta_k + dd^c u_k)^n = c_k (1 + f) \omega^n.
			\]
			Integrating both sides, we obtain that
			\[
			c_k \geq \frac{\underline{\text{Vol}}(\beta_k)}{\int_X (1 + f) \omega^n} \geq \frac{\underline{\text{Vol}}(\beta)}{\int_X (1 + f) \omega^n} > 0.
			\]
			On the other hand, for any $k \geq j$,
			\[
			(\beta_j + dd^c u_k)^n \geq (\beta_k + dd^c u_k)^n \geq \frac{\underline{\text{Vol}}(\beta)}{\int_X (1 + f) \omega^n} f \omega^n.
			\]
			By the maximum principle (\cref{max principle 3}), it follows that
			\[
			(\beta_j + dd^c \Phi_j)^n \geq \frac{\underline{\text{Vol}}(\beta)}{\int_X (1 + f) \omega^n} f \omega^n.
			\]
			The result follows by letting $j \rightarrow \infty$.
		\end{proof}   
		\section{Solving complex Monge-Amp\`ere equations}

		Following \cite{GL23}, we can solve the complex  Monge-Amp\`ere equation with an exponential twist:
		\begin{theorem}\label{MA equation 1}
			Let $(X,\omega)$ be a compact Hermitian manifold of complex dimension $n$, equipped with a Hermitian metric $\omega$. 
			Let $\beta$ be  a possibly non-closed smooth $(1,1)$-form such that there exists a bounded $\beta$-psh function $\rho$ and $\underline{\text{Vol}}(\beta) > 0$.
			Assume $0\leq f\in L^p(X,\omega^n)$ with $p>1$ and $||f||_p>0$. Then, for each $\lambda>0$, there exists a unique $\varphi_{\lambda}\in \mbox{PSH}(X,\beta)\cap L^{\infty}(X)$ such that
			$$
			(\beta+dd^c\varphi_{\lambda})^n=e^{\lambda\varphi_{\lambda}}f\omega^n.
			$$
			Furthermore, we have $|\varphi_{\lambda}|\leq C$ for some uniform constant C depending only on $\lambda,\beta,p,||f||_p$, $X,\omega$.
		\end{theorem}
		
		\begin{proof}
			We can assume without loss of generality that $\lambda=1$. Set again $\beta_j:=\beta+\varepsilon_j\omega$ for some decreasing sequence of  $\varepsilon_j\searrow0$. Using \cite[Theorem 0.1]{Ngu16} (see also \cite[Theorem 4.5]{BGL24})  we can solve
			$$
			(\beta_j+dd^c\varphi_j)^n=e^{\varphi_j}f\omega^n
			$$
			for some $\varphi_j\in \mbox{PSH}(X,\beta_j)\cap L^{\infty}(X)$. By assumption $\beta_j>\beta_{j+1}$, and we have
			$$
			e^{-\varphi_j}(\beta_j+dd^c\varphi_j)^n=f\omega^n=e^{-\varphi_{j+1}}(\beta_{j+1}+dd^c\varphi_{j+1})^n\leq e^{-\varphi_{j+1}}(\beta_{j}+dd^c\varphi_{j+1})^n.
			$$
			By the domination principle (\cref{domination for beta 2}), we have $\varphi_{j+1}\leq\varphi_j$ and hence the sequence $\varphi_j$ decreases to some $\varphi$.
			
			We claim that $\varphi_j$ is uniformly bounded and hence $\varphi\in \mbox{PSH}(X,\beta)$. By \cref{subsolution}, we can find a subsolution $u$ satisfying
			$$
			(\beta+dd^cu)^n\geq cf\omega^n\geq e^{u+\ln c}f\omega^n
			$$
			since $u\leq0$. This implies that
			$$
			e^{-(u+\ln c)}(\beta+dd^c(u+\ln c))^n\geq f\omega^n=e^{-\varphi_j}(\beta_j+dd^c\varphi_j)^n
			$$
			and thus $\varphi_j\geq u+\ln c$ by the domination principle  (\cref{domination for beta 2}). The claim follows.
			
			By the Bedford-Taylor convergence theorem it follows  that $(\beta+dd^c\varphi)^n=e^{\varphi}f\omega^n$. The uniqueness is again a consequence of the domination principle (\cref{domination for beta 2}).
		\end{proof}

		\begin{theorem}\label{main thm}
			Let $(X,\omega)$ be a compact Hermitian manifold of complex dimension $n$, equipped with a Hermitian metric $\omega$. 
			Let $\beta$ be  a possibly non-closed smooth $(1,1)$-form such that there exists a bounded $\beta$-psh function $\rho$ and $\underline{\text{Vol}}(\beta) > 0$.
			Assume $0\leq f\in L^p(X,\omega^n)$ with $p>1$ and $||f||_p>0$. Then there exists a constant $c>0$ and a function $\varphi\in \mbox{PSH}(X,\beta)\cap L^{\infty}(X)$ such that
			$$
			(\beta+dd^c\varphi)^n=cf\omega^n.
			$$
			Furthermore, the constant $c$ is uniquely determined by $f,\beta,X$ and $\mbox{osc} \varphi\leq C$ for some constant $C$ depending only on $\beta,p,||f||_{p},X,\omega$.
		\end{theorem}
		\begin{proof}
			
			The proof follows the argument presented in \cite[Theorem 4.6]{BGL24}.
			
			The uniqueness of the constant $c$ is an immediate consequence of the domination principle (\cref{domination for beta 3}). For the existence, we first apply \cref{MA equation 1} to obtain $u_j \in \text{PSH}(X, \beta) \cap L^{\infty}(X)$ satisfying
			\[
			(\beta + dd^c u_j)^n = e^{j^{-1} u_j} f \omega^n.
			\]
			Define $v_j := u_j - \sup_X u_j$ and $f_j := e^{j^{-1} v_j} f$. The above equation then transforms into
			\[
			(\beta + dd^c v_j)^n = c_j f_j \omega^n,
			\]
			where $c_j := e^{j^{-1} \sup_X u_j}$. By extracting a subsequence if necessary, we may assume that $v_j \rightarrow v \in \text{PSH}(X, \beta)$ in $L^1(X,\omega^n)$ and almost everywhere, implying that $e^{j^{-1} v_j} \rightarrow 1$ in $L^1(X,\omega^n)$.
			
			To establish a uniform upper bound for $c_j$, we select a Gauduchon metric $\omega_G$ on $X$ and write $\omega_G^n = h \omega^n$ for some positive function $h\in C^\infty(X)$. After rescaling, we may assume $\beta \leq \omega_G$, leading to $(\omega_G + dd^c v_j)^n \geq (\beta + dd^c v_j)^n \geq c_j f_j \omega^n$. Applying the mixed-type inequality (\cref{mixed type 1}), we deduce
			\[
			(\omega_G + dd^c v_j) \wedge \omega_G^{n-1} \geq c_j^{\frac{1}{n}} f^{\frac{1}{n}} h^{1 - \frac{1}{n}} \omega^n.
			\]
			Integrating both sides as in the proof of \cref{subsolution}, we obtain a uniform upper bound for $c_j$.
			
			By Skoda's uniform integrability theorem (\cite[Theorem 8.11]{GZ17}), we choose a small constant $\varepsilon > 0$ such that $e^{-\varepsilon v_j} f$ is uniformly in $L^q$ for some $1 < q < p$. Using \cref{subsolution}, we find a subsolution $\eta_j$ satisfying
			\[
			(\beta + dd^c \eta_j)^n \geq m e^{-\varepsilon v_j} f \omega^n \geq m e^{(\varepsilon + \frac{1}{j}) \eta_j - \varepsilon v_j} f \omega^n.
			\]
			From the first inequality, we have
			\[
			(\beta + dd^c \eta_j)^n \geq m e^{-\varepsilon v_j} f \omega^n = (\beta + dd^c v_j)^n \frac{m}{c_j} e^{-(\varepsilon + j^{-1}) v_j}.
			\]
			Since $v_j \leq 0$, the domination principle (\cref{domination for beta 3}) implies $m \leq c_j$.
			
			From the second inequality, we obtain
			\[
			e^{-(\varepsilon + \frac{1}{j}) \eta_j} (\beta + dd^c \eta_j)^n \geq e^{-(\varepsilon + \frac{1}{j}) (v_j - \frac{\ln \frac{m}{c_j}}{\varepsilon + \frac{1}{j}})} (\beta + dd^c v_j)^n.
			\]
			Applying the domination principle (\cref{domination for beta 2}), we get $v_j \geq \eta_j + \frac{\ln \frac{m}{c_j}}{\varepsilon + \frac{1}{j}} \geq \eta_j - C \geq -T - C$ by \cref{subsolution} and the fact that $c_j$ is uniformly bounded above. This shows that $v_j$ is uniformly bounded below. Moreover, since $c_j$ has uniform upper and lower bounds, we may assume, by extracting a subsequence, that $c_j \rightarrow c > 0$.
			
			Define
			\[
			\Phi_j := P_{\beta} \left( \inf_{k \geq j} v_k \right), \quad \psi_j := \left( \sup_{k \geq j} v_k \right)^*.
			\]
			From the preceding discussion, $\Phi_j$ and $\psi_j$ are uniformly bounded, and $\psi_j \searrow v$ by Hartogs' lemma. Using a standard argument involving \cref{max principle 3} for finitely many functions and then taking the limit, we obtain
			\[
			(\beta + dd^c \Phi_j)^n \leq e^{\varepsilon (\Phi_j - \inf_{k \geq j} v_k)} \sup_{k \geq j} (c_k f_k) \omega^n
			\]
			and
			\[
			(\beta + dd^c \psi_j)^n \geq e^{\varepsilon (\psi_j - \sup_{k \geq j} v_k)} \inf_{k \geq j} (c_k f_k) \omega^n.
			\]
			Suppose $\Phi_j$ increases to some $\Phi \in \text{PSH}(X, \beta)$. Clearly, $\Phi \leq v$. Letting $j \rightarrow \infty$, we get
			\[
			(\beta + dd^c \Phi)^n \leq e^{\varepsilon (\Phi - v)} c f \omega^n \quad \text{and} \quad (\beta + dd^c v)^n \geq c f \omega^n \quad (\text{since } e^{\varepsilon (v - v)} = 1).
			\]
			Applying the domination principle once more, we conclude $\Phi \geq v$ and hence $\Phi = v$, which implies $(\beta + dd^c v)^n = c f \omega^n$.
		\end{proof}
		
		Following the method of \cite{GL23, LPT20}, (see also \cite{GL25} for a version of complex Hessian equations) we establish a stability result of the equation in \cref{MA equation 1}.

		\begin{theorem}\label{thm:stablity of MA 1}
			Let $(X, \omega)$ be a compact Hermitian manifold of complex dimension $n$, equipped with a Hermitian metric $\omega$. Let $\beta$ be a possibly non-closed smooth $(1,1)$-form on $X$ such that there exists a bounded  $\beta$-psh function $\rho$ and $\underline{\mathrm{Vol}}(\beta) > 0$. 
			Fix $f, g \in L^p(X, \omega^n)$ with $p > 1$ and $\|f\|_p, \|g\|_p>0$. Let $\varphi, \psi \in \mathrm{PSH}(X, \beta) \cap L^{\infty}(X)$ be such that
			\[
			(\beta + dd^c \varphi)^n = e^{\lambda \varphi} f \omega^n, \quad \text{and} \quad (\beta + dd^c \psi)^n = e^{\lambda \psi} g \omega^n,
			\]
			where $\lambda>0$ is a fixed constant. 
			Then we have the estimate
			\[
			\|\varphi - \psi\|_{\infty} \leq C \|f - g\|_p^{\frac{1}{n}},
			\]
			where $C > 0$ is a constant depending on $n$, $p$, $\lambda$, $\|f\|_p$, and $\|g\|_p$. 
			In particular, this implies the uniqueness of the solution to the Monge-Amp\`ere equation in \cref{MA equation 1}.
		\end{theorem}
		\begin{proof}
			Up to rescaling, we may assume without loss of generality that $\lambda = 1$. Firstly, if $\|f - g\|_p = 0$, by the uniqueness property, we deduce that $\varphi = \psi$, and there is nothing to prove. Thus, we may assume $\|f - g\|_p > 0$. Applying \cref{main thm}, we can solve
			\[
			(\beta + dd^c u)^n = c\left(\frac{|f - g|}{\|f - g\|_p} + 1\right)\omega^n := ch\omega^n
			\]
			with $\sup_X u = 0$. Integrating both sides of the equation, we obtain
			\[
			c\|h\|_1 \geq \underline{\text{Vol}}(\beta),
			\]
			and thus $c \geq c_1 > 0$ is uniformly bounded below.
			
			Pick a Gauduchon metric $\omega_G$ on $X$. We can assume that $\omega_G \geq \beta$. It follows that $u \in \text{PSH}(X, \omega_G)$ and $(\omega_G + dd^c u)^n \geq ch\delta\omega_G^n$ for some uniform constant $\delta > 0$. Since $h \geq 1$, the mixed-type inequality (\cref{mixed type 1}) yields that
			\[
			(\omega_G + dd^c u) \wedge \omega_G^{n-1} \geq c^{\frac{1}{n}}\delta^{\frac{1}{n}}\omega_G^n.
			\]
			Consequently, we obtain an upper bound for $c$.
			
			Now, we divide the proof into two cases: Set $\varepsilon := e^{\frac{\sup_X \varphi - \log c}{n}}\|f - g\|_p^{\frac{1}{n}}$. If $\varepsilon \geq\frac{1}{2}$, then $\|f - g\|_p^{\frac{1}{n}} \geq \frac{c_1^{\frac{1}{n}}}{2}e^{\frac{-\sup_X \varphi}{n}}$. Since $\varphi, \psi$ are uniformly bounded by a constant $C_1$ by \cref{MA equation 1}, we just need to choose $C \geq 4C_1 c^{-\frac{1}{n}}e^{\frac{\sup_X \varphi}{n}}$.
			
			Hence, we can assume that $\varepsilon < \frac{1}{2}$. Set
			\[
			\phi := (1 - \varepsilon)\varphi + \varepsilon u - K\varepsilon + n\log(1 - \varepsilon)
			\]
			with a constant $K$ to be defined later. Then we compute:
			
			\begin{align*}
				(\beta + dd^c \phi)^n &\geq (1 - \varepsilon)^n e^{\varphi} f \omega^n + \varepsilon^n ch \omega^n \\
				&= \left[(1 - \varepsilon)^n e^{\varphi} f + e^{\sup_X \varphi - \log c} \|f - g\|_p \, c \left(\frac{|f - g|}{\|f - g\|_p} + 1\right)\right] \omega^n \\
				&\geq \left[(1 - \varepsilon)^n e^{\varphi} f + e^{\sup_X \varphi} |f - g|\right] \omega^n \\
				&\geq e^{\varphi + n\log(1 - \varepsilon)} f \omega^n + e^{\varphi} |f - g| \omega^n \\
				&\geq e^{\varphi + n\log(1 - \varepsilon)} g \omega^n.
			\end{align*}
			
			Choose $K := \sup_X (-\varphi) \geq u - \varphi$. We easily obtain that
			\[
			(\beta + dd^c \phi)^n \geq e^{\phi} g \omega^n.
			\]
			It thus follows from the domination principle (\cref{domination for beta 2}) that $\phi \leq \psi$ and hence $\varphi \leq \psi + C_2 \varepsilon$ for a uniform constant $C_2$. Reversing the inequality, we conclude the proof.

			Up to rescaling we may assume without loss of generality that $\lambda=1$. Firstly, if $||f-g||_p=0$, by the uniqueness we deduce that $\varphi=\psi$ and there is nothing to prove. So we may assume $||f-g||_p>0$. Applying \cref{main thm} we can solve
			$$
			(\beta+dd^cu)^n=c\left(\frac{|f-g|}{||f-g||_p} +1\right)\omega^n:=ch\omega^n 
			$$
			with $\sup_Xu=0$. Integrating both sides of the equation, we have
			$$
			c||h||_1\geq\underline{\mbox{Vol}}(\beta)
			$$
			and thus $c\geq c_1>0$ is uniformly bounded below.
			
			Pick a Gauduchon metric $\omega_G$ on $X$, we can assume that $\omega_G\geq\beta$. It follows that $u\in \mbox{PSH}(X,\omega_G)$ and $(\omega_G+dd^cu)^n\geq ch\delta\omega_G^n$ for some uniform constant $\delta>0$. Since $h\geq1$, the mixed type inequality \cref{mixed type 1} yields that
			$$
			(\omega_G+dd^cu)\wedge\omega_G^{n-1}\geq c^{\frac{1}{n}}\delta^{\frac{1}{n}}\omega_G^n.
			$$
			Consequently, we get an upper bound for $c$.
			
			Now we divide the proof into two cases: set $\varepsilon:=e^{\frac{\sup_X\varphi-\log\,c}{n}}||f-g||_p^{\frac{1}{n}}$. If $\varepsilon>\frac{1}{2}$, then $||f-g||_p^{\frac{1}{n}}\geq\frac{c_1^{\frac{1}{n}}}{2}e^{\frac{-\sup_X\varphi}{n}}$. Since $\varphi,\psi$ are uniformly bounded by a constant $C_1$ by \cref{MA equation 1}, we just need to choose $C\geq4C_1c^{-\frac{1}{n}}e^{\frac{\sup_X\varphi}{n}}$.
			
			Hence we can assume that $\varepsilon<\frac{1}{2}$. Set 
			$$
			\phi:=(1-\varepsilon)\varphi+\varepsilon u-K\varepsilon+n\log(1-\varepsilon)
			$$
			with a constant $K$ to be defined later.
			Then we compute
			
			\begin{align*}
				(\beta+dd^c\phi)^n&\geq(1-\varepsilon)^ne^{\varphi}f\omega^n+\varepsilon^nch\omega^n\\
				&=\left[(1-\varepsilon)^ne^{\varphi}f+e^{\sup_X\varphi-\log\,c}||f-g||_p\,c\left(\frac{|f-g|}{||f-g||_p} +1\right)\right]\omega^n\\
				&\geq\left[(1-\varepsilon)^ne^{\varphi}f+e^{\sup_X\varphi}|f-g|\right]\omega^n\\ 
				&\geq e^{\varphi+n\log(1-\varepsilon)}f\omega^n+e^{\varphi}|f-g|\omega^n\\
				&\geq e^{\varphi+n\log(1-\varepsilon)}g\omega^n.
			\end{align*}
			Choose $K:=\sup_X(-\varphi)\geq u-\varphi$, we obtain easily that
			$$
			(\beta+dd^c\phi)^n\geq e^{\phi}g\omega^n.     
			$$     
			It thus follows from the domination principle (\cref{domination for beta 2}) that $\phi\leq\psi$ and hence $\varphi\leq\psi+C_2\varepsilon$ for a uniform constant $C_2$. Reversing the inequality we conclude the proof.
		\end{proof}
		
		\begin{remark}
			The uniqueness of solutions in \cref{main thm} is a subtle issue. The methods in \cite{KN19}, \cite{LPT20} and \cite{GL23} do not seem directly applicable to our setting, even in the case where $f\geq c_0>0$. Interested readers may attempt to adapt the approaches from \cite[Corollary 3.9]{KN19} and \cite[Theorem 3.5]{GL23} to address this problem.
		\end{remark}

		\section{Asymptotics  of solutions of  degenerate complex Monge-Amp\`ere equations}
		Following the method of \cite[Theorem 2.6]{GL21}, we wish to establish the following stability result. This result will yield a $L^1-L^{\infty}$ stability estimate which will be useful in the proof of the Demailly-P\v{a}un's conjecture (see \cref{thm:DP conjecture}):
		\begin{theorem}\label{thm:stability}
			Let $(X,\omega)$ be a compact Hermitian manifold of complex dimension $n$, equipped with a Hermitian metric $\omega$. 
			Let $\beta$ be  a possibly non-closed smooth $(1,1)$-form such that there exists a bounded $\beta$-psh function $\rho$ and $\underline{\text{Vol}}(\beta) > 0$. Let $\mu$ be a probability measure on X such that $\mbox{PSH}(X,\beta)\subset L^m(\mu)$ for some $m>n$. Let $\varphi\in \mbox{PSH}(X,\beta)\cap L^{\infty}(X)$ be such that $\sup_X\varphi=0$ and $(\beta+dd^c\varphi)^n\leq \tau\mu$ for some $\tau>0$. Then
			$$
			\underset{X}{\sup}(\phi-\varphi)_+\leq T\left(\int_X(\phi-\varphi)_+d\mu\right)^{\gamma}    
			$$
			for any $\phi\in \mbox{PSH}(X,\beta)\cap L^{\infty}(X)$. Where $\gamma=\gamma(m,n)>0$ and $T=T(\mu,||\phi||_{\infty})$ is a uniform constant which depends on an upper bound of $||\phi||_{\infty},\frac{\tau}{\underline{\mbox{Vol}}(\beta)}$ and on
			$$
			A_m(\mu):=\sup\left\{\left(\int_X(-\psi)^md\mu\right)^{\frac{1}{m}}: \psi\in \mbox{PSH}(X,\beta)\,with\,\underset{X}{\sup}\,\psi=0\right\}.
			$$
		\end{theorem}
		\begin{proof}
			We can assume without loss of generality that $\varphi\leq\phi$. Let $T_{\max}:=\sup\{t>0:\mu\{\varphi<\phi-t\}>0\}$. By \cref{thm: a prior estimate} we have $T_{\max}$ is uniformly bounded by $\mu,||\phi||_{\infty}$ and $\frac{c}{\underline{\mbox{Vol}}(\beta)}$. Let $\chi:\mathbb{R}^-\rightarrow\mathbb{R}^-$ be a concave increasing weight such that $\chi(0)=0$ and $\chi'(0)=1$. Set $\psi:=\chi\circ(\varphi-\phi)+\phi$ and $u:=P_{\beta}(\psi)$. By \cref{lem: concave weight} we have
			$$
			\frac{1}{\tau}(\beta+dd^cu)^n\leq\mathds{1}_{\{u=\psi\}}(\chi'\circ(\varphi-\phi))^n\mu.
			$$
			\textbf{Step 1.} We first control the energy of $u$. Fix $0<a<b<c<2c<\varepsilon$ so small that $q:=\frac{(\varepsilon-a)(n+b)}{b-a}<m=n+\varepsilon$. It follows that $\mbox{PSH}(X,\beta)\subset L^{n+2c}(\mu)$. Since $\chi(0)=0$ and $\chi$ is concave, we have $|\chi(t)|\leq |t|\chi'(t)$. We can thus estimate the energy :
			\begin{align*}
				0\leq\int_X(\phi-u)^{c}\frac{(\beta+dd^cu)^n}{\tau}\leq&\int_X(-\chi\circ(\varphi-\phi))^{c}(\chi'\circ(\varphi-\phi))^nd\mu\leq\int_X(\phi-\varphi)^{c}(\chi'\circ(\varphi-\phi))^{n+c}d\mu\\    
				\leq&\left(\int_X(\phi-\varphi)^{n+2c}d\mu\right)^{\frac{c}{n+2c}}\left(\int_X(\chi'\circ(\varphi-\phi))^{n+2c}d\mu\right)^{\frac{n+c}{n+2c}}\\
				\leq&\left(\int_X(||\phi||_{\infty}-\varphi)^{n+2c}d\mu\right)^{\frac{c}{n+2c}}\left(\int_X(\chi'\circ(\varphi-\phi))^{n+2c}d\mu\right)^{\frac{n+c}{n+2c}}\\
				\leq&(A_m(\mu)+||\phi||_{\infty})^{c}\left(\int_X(\chi'\circ(\varphi-\phi))^{n+2c}d\mu\right)^{\frac{n+c}{n+2c}}.
			\end{align*}
			Where the first inequality follows since $u=\chi\circ(\varphi-\phi)+\phi$ on the support of $(\beta+dd^cu)^n$, the second inequality follows from the fact $|\chi(t)|\leq |t|\chi'(t)$, the third follows from H\"older's inequality while in the last we used the fact that $\mu$ is a probability measure and the monotonicity of $L^p$-norms with respect to $\mu$.
			
			\textbf{Step 2.} We choose an appropriate weight $\chi$ such that $\int_X(\chi'\circ(\varphi-\phi))^{n+2c}d\mu\leq B$ is under control.  We first fix a $T\in[0,T_{max})$. Recall that if $g:\mathbb{R}^+\rightarrow\mathbb{R}^+$ is an increasing, absolutely continuous function (on $[0,T]$) with $g(0)=1$, then
			
			\[\int_X g\circ(\phi-\varphi)d\mu=\mu(X)+{\int_0^{T_{\max}}} g'(t)\mu(\varphi<\phi-t)dt.\]
			Set $f(t):=\begin{cases}  
				\frac{1}{\mu(\varphi<\phi-t)}, & \text{if } t \in [0,T] \\
				1, & \text{if } t >T , 
			\end{cases}$ , then $f(t)$ is a $L^1$ function. Set $g(x):=\int_0^xf(t)dt+1$. It is easy to check that g is absolutely continuous on $[0,T]$ and $g'(t)=f(t)$ almost everywhere ( at least at every Lebesgue point). 
			
			Using the argument again, by letting $-\chi(-x):=\int_0^xg(t)^{\frac{1}{n+2c}}dt$, then $\chi(0)=0,\chi'(0)=1$ and $g(t)=(\chi'(-t))^{n+2c}$, note that this equality holds everywhere since $g$ is continuous and $g\geq1$. The choice guarantees that $\chi$ is concave increasing with $\chi'\geq1$, and
			$$
			\int_X(\chi'\circ(\varphi-\phi))^{n+2c}d\mu\leq \mu(X)+\int_0^{T_{\max}}1dt=1+T_{\max}.
			$$
			Since $T_{\max}$ is uniformly bounded above by \cref{thm: a prior estimate}, we may choose $B:=1+T_{\max}$.
			
			\textbf{Step 3.} We next turn to estimate the level set $\mu(\varphi<\phi-t)$ by controlling the norm of $||u||_{L^m}$. Since $g\geq1$, our construction in \textbf{Step 2} yields that $\chi(x)\leq x$ and hence $u\leq\phi+\chi(\varphi-\phi)\leq\varphi\leq0$.  Moreover,
			\begin{align*}
				0\leq(-\underset{X}{\sup}(u-\phi))^{c}\frac{\underline{\mbox{Vol}}(\beta)}{\tau}\leq&(-\underset{X}{\sup}(u-\phi))^{c}\int_X\frac{(\beta+dd^cu)^n}{\tau}\\
				\leq&\int_X(\phi-u)^{c}\frac{(\beta+dd^cu)^n}{\tau}\\
				\leq&(A_m(\mu)+||\phi||_{\infty})^{c}B^{\frac{n+c}{n+2c}}.
			\end{align*}
			This gives a uniform lower bound $C_1$ for $\underset{X}{\sup}(u-\phi)$. We furthermore have $0\geq \sup_Xu\geq \sup_X(u-\phi)+\inf_X\phi\geq C_1-||\phi||_{\infty}$ and hence $u$ lies  in a compact set of $\mbox{PSH}(X,\beta)$. We can write
			\begin{align*}
				\int_X|\chi(\varphi-\phi)|^md\mu&\leq\int_X|\chi(\varphi-\phi)|^{n+a}(\phi-u)^{\varepsilon-a}d\mu\\
				&\leq\left(\int_X|\chi(\varphi-\phi)|^{n+b}d\mu\right)^{\frac{n+a}{n+b}}\left(\int_X(\phi-u)^qd\mu\right)^{\frac{b-a}{n+b}}\\
				&\leq C_2\left(\int_X|(\phi-\varphi)\chi^{\prime}(\varphi-\phi)|^{n+b}d\mu\right)^{\frac{n+a}{n+b}}\\
				&\leq C_2\left(\int_X|\chi^{\prime}(\varphi-\phi)|^{n+c}d\mu\right)^{\frac{n+a}{n+c}}\left(\int_X(\phi-\varphi)^{\frac{n+c}{c-b}}d\mu\right)^{\gamma}\\
				&\leq C_3B^{\frac{n+a}{n+2c}}\left(\int_X(\phi-\varphi)d\mu\right)^{\gamma}:=C_4\delta,
			\end{align*}
			where $\gamma:=\frac{(c-b)(n+a)}{(n+b)(n+c)}$ and $\delta:=\left(\int_X(\phi-\varphi)d\mu\right)^{\gamma}$. The first inequality is because $u-\phi\leq\chi(\varphi-\phi)\leq0$ and $m:=n+\varepsilon$, the second and the forth inequality follows from the Hölder's inequality and the third follows from the fact $|\chi(t)|\leq |t|\chi^{\prime}(t)$.
			
			Therefore, we have
			\begin{align*}
				\mu(\varphi<\phi-t)=&\int_{\{\varphi<\phi-t\}}d\mu\leq\int_X\frac{|\chi\circ(\varphi-\phi)|^m}{|\chi(-t)|^m}d\mu\\
				=&\frac{1}{|\chi(-t)|^m}||\chi\circ(\varphi-\rho)||_{L^m(\mu)}^m\leq \frac{C_4\delta}{|\chi(-t)|^m}.
			\end{align*}
			\textbf{Step 4.} Conclusion. Set $h(t):=-\chi(-t)$. We have $h(0)=0 $ and $h'(t)=[g(t)]^{\frac{1}{n+2c}}$ is positive increasing, hence $h$ is convex. Observe also that $g(t)\geq g(0)\geq1$, we have $h'(t)\geq1$ and $h(1)=\int_0^1h'(t)dt\geq 1$. For almost all $t\in[0,T]$,
			$$
			\frac{1}{g'(t)}=\mu(\varphi<\phi-t)\leq \frac{C_4\delta}{h^m(t)}.
			$$
			This implies that 
			$$
			h^m(t)\leq C_4\delta g'(t)=(n+2c)C_4\delta h''(t)(h')^{n+2c-1}(t)
			$$
			almost everywhere.
			It follows that we can multiply both sides by $h'$ and integrate between $0$ and $t$, we obtain
			\begin{align*}
				\frac{h^{m+1}(t)}{m+1}\leq&(n+2c)C_4\delta\int_0^th''(s)(h')^{n+2c}(s)ds\\
				\leq&C_5\delta((h')^{n+2c+1}(t)-1).
			\end{align*}
			Note that the first and the second inequality holds because $h$ and $h'$ are absolutely continuous and hence we can take derivative and integrate. This yields that
			$$
			1\leq\frac{C_6\delta(h')^{n+2c+1}(t)}{C_6\delta+h(t)^{m+1}}.
			$$
			Set $\alpha:=m+1=n+\varepsilon+1>\eta:=n+2c+1>2$.
			The previous inequality then reads $1\leq\frac{C_7\delta^{\frac{1}{\eta}} h'(t)}{(C_6\delta+h(t)^{\alpha})^{\frac{1}{\eta}}}$. Integrating both sides between $0$ and $T$ we get
			$$
			T\leq C_7\delta^{\frac{1}{\eta}}\int_0^T\frac{ h'(t)}{(C_6\delta+h(t)^{\alpha})^{\frac{1}{\eta}}}dt=C_7\int_0^{h(T)\delta^{-\frac{1}{\alpha}}}\frac{\delta^{\frac{1}{\alpha}}dx}{(C_6+x^{\alpha})^{\frac{1}{\eta}}}\leq C_7\delta^{\frac{1}{\alpha}}\int_0^{+\infty}\frac{dx}{(C_6+x^{\alpha})^{\frac{1}{\eta}}}\leq C_8\delta^{\frac{1}{\alpha}},
			$$
			where in the second equality we have used the change of variable $x=h(t)\delta^{-\frac{1}{\alpha}}$.
			Hence
			$$
			T\leq C_8\delta^{\frac{1}{\alpha}}\leq C_8(\int_X(\phi-\varphi)_+d\mu)^{\frac{\gamma}{\alpha}}
			$$
			by the definition of $\delta$.
			Finally, letting $T\rightarrow T_{\max}$ we conclude the proof.

		\end{proof}
		
		\begin{corollary}\label{cor:stability 2}
			Let $\beta,\mu_1,\mu_2$ be as in \cref{thm: a prior estimate}. Let $\varphi,\phi\in \mbox{PSH}(X,\beta)\cap L^{\infty}(X)$ be such that $\sup_X\varphi=0,\sup
			_X\phi=0$ and $(\beta+dd^c\varphi)^n\leq \tau_1\mu_1,(\beta+dd^c\phi)^n\leq \tau_2\mu_2$ for some $\tau_1>0,\tau_2>0$. Then
			$$
			||\varphi-\phi||_{\infty}\leq T||\varphi-\phi||_{L^1(\mu_1+\mu_2)}^{\gamma},
			$$
			where $\gamma=\gamma(m,n)>0$ is a constant and $T$ is a uniform constant depending on the upper bound of $||\phi||_{\infty},||\varphi||_{\infty},\frac{\max(\tau_1,\tau_2)}{\underline{\mbox{Vol}}(\beta)}$ and $A_m(\mu_1),A_m(\mu_2)$.
		\end{corollary}
		\begin{proof}
			This is a direct consequence of \cref{thm:stability}.
		\end{proof}
		\section{Applications}
		\subsection{The extended Tosatti-Weinkove conjecture}
		
		In this section, we establish an extended version of the Tosatti-Weinkove conjecture \cite[Conjectures 1.1] {LWZ24} when the background cohomology class $\{\beta\}$ is non-closed, extending the result in \cite[Theorem 1.10]{LWZ24}:
		
		\begin{theorem}\label{thm: TW conj}
			Let $X$ be a compact Hermitian manifold. Let $\beta\in BC^{1,1}(X)$ be as in the introduction. Let $x_1,...,x_N$ be fixed points and $\tau_1,...,\tau_N$ be positive real numbers such that
			$$
			\sum_{i=1}^{N}\tau_i^n<\underline{\mbox{Vol}}(\beta).    
			$$
			Then there exists a $\beta$-psh function $\varphi$ with logarithmic poles at $x_1,...,x_N$:
			$$
			\varphi(z)\leq c\tau_j\log|z|+O(1)
			$$
			in a coordinate ball $(z_1,...,z_n)$ centered at $x_j$, where $c\geq1$ is a uniform constant.
		\end{theorem}
		\begin{proof}
			The argument mainly follows the idea from \cite{LWZ24}. Let $\chi:\mathbb{R}\rightarrow\mathbb{R}$ be a smooth convex increasing function such that $\chi(t)=t$ for $t\geq0$ and $\chi(t)=-\frac{1}{2}$ for $t\leq-1$. For a fixed $j$, let $(z_1,...,z_n)$ be a coordinate chart centered at $x_j$. We define
			$$
			\gamma_{j,\varepsilon}:=dd^c(\chi\circ(\log\frac{|z|}{\varepsilon}))
			$$
			for any $\varepsilon>0$. By the choice of $\chi$, it is clear that $\gamma_{j,\varepsilon}$ is a smooth closed positive $(1,1)$ form on this coordinate chart. Observe also that $\gamma_{j,\varepsilon}=dd^c\log|z|$ outside the set $\{|z|<\varepsilon\}$ and hence $\gamma_{j,\varepsilon}^n=0$ outside $\{|z|<\varepsilon\}$. Since
			$$
			\int_X\gamma_{j,\varepsilon}^n=1,
			$$
			we have $\gamma_{j,\varepsilon}^n\rightarrow\delta_{x_j}$ as $\varepsilon\rightarrow0$. Set
			$$
			\delta:=\underline{\mbox{Vol}}(\beta)-\sum_{i=1}^{N}\tau_i^n>0.
			$$
			By \cref{main thm}, we can solve the following family of degenerate CMA equation
			\begin{equation}\label{eq:solve eq}
				(\beta+dd^c\varphi_{\varepsilon})^n=c_{\varepsilon}\left(\sum_{i=1}^{N}\tau_i^n\gamma_{j,\varepsilon}^n+\delta\frac{\omega^n}{\int_X\omega^n}\right) 
			\end{equation}
			for some constant $c_{\varepsilon}>0$ with $\varphi_{\varepsilon}\in \mbox{PSH}(X,\beta)\cap L^{\infty}(X)$ and $\sup_X\varphi_{\varepsilon}=0$. We first claim that the family $c_{\varepsilon}$ is uniformly bounded.
			
			Integrate both sides of \cref{eq:solve eq} we obtain that
			$$
			c_{\varepsilon}\underline{\mbox{Vol}}(\beta)=\int_X(\beta+dd^c\varphi_{\varepsilon})^n\geq\underline{\mbox{Vol}}(\beta)
			$$
			and hence $c_{\varepsilon}\geq1$. For the upper bound, pick a Gauduchon metric $\omega_G$ on $X$ such that $\omega_G\geq\beta$. We then have $(\omega_G+dd^c\varphi_{\varepsilon})^n\geq c_{\varepsilon}\delta^{\prime}\omega_G^n$ for some positive constant $\delta^{\prime}>0$. Using the mixed type inequality (\cref{mixed type 1}) we have
			$$
			(\omega_G+dd^c\varphi_{\varepsilon})\wedge\omega_G^{n-1}\geq c_{\varepsilon}^{\frac{1}{n}}\delta^{\prime\frac{1}{n}}\omega_G^n.
			$$
			Integrating both sides, we clearly get an upper bound for $c_{\varepsilon}$.
			
			Since we have assumed that $\sup_X\varphi_{\varepsilon}=0$, up to extracting a subsequence we can assume that $\varphi_{\varepsilon}$ converges to a function $\varphi\in \mbox{PSH}(X,\beta)$ in $L^{1}(X,\omega^n)$ and almost everywhere, and $c_{\varepsilon}\rightarrow c\geq1$. We claim that $\varphi$ has the desired singularities.
			
			Let $U$ be a neighborhood of $x_j$ and choose a smooth function $h$ near $U$ such that $\beta\leq dd^ch$. Set $v_{\varepsilon}:=h+\varphi_{\varepsilon}\in \mbox{PSH}(U)$. Since we have assumed that $\sup_X\varphi_{\varepsilon}=0$, there is a uniform constant $C>0$ such that $v_{\varepsilon}|_{\partial U}\leq C$. Set 
			$$
			u_{\varepsilon}:=c_{\varepsilon}^{\frac{1}{n}}\tau_j(\chi\circ(\log\frac{|z|}{\varepsilon}))+\log\varepsilon)+C_1
			$$
			for a large constant $C_1$.
			Then by the choice of $\chi$, for $\varepsilon$ small enough we have
			$$
			u_{\varepsilon}|_{\partial U}=c_{\varepsilon}^{\frac{1}{n}}\tau_j\log|z|+C_1,\quad v_{\varepsilon}|_{\partial U}\leq C,
			$$
			and
			$$
			(dd^cv_{\varepsilon})^n\geq(\beta+dd^c\varphi_{\varepsilon})^n=c_{\varepsilon}\left(\sum_{i=1}^{N}\tau_i^n\gamma_{j,\varepsilon}^n+\delta\frac{\omega^n}{\int_X\omega^n}\right)>c_{\varepsilon}\sum_{i=1}^{N}\tau_i^n\gamma_{j,\varepsilon}^n\geq(dd^cu_{\varepsilon})^n.
			$$
			Since we have proved that $c_\varepsilon$ is uniformly bounded (with respect to  $\varepsilon$), we can choose $C_1$ large enough such that $u_{\varepsilon}>v_{\varepsilon}$ on $\partial U$ for all $\varepsilon$. Then by the Bedford-Taylor comparison principle \cite[Theorem 4.1]{BT82}, we get that $u_{\varepsilon}\geq v_{\varepsilon}$ on $U$. Thus we get
			$$
			\varphi_{\varepsilon}\leq c_{\varepsilon}^{\frac{1}{n}}\tau_j\left(\chi\circ(\log\frac{|z|}{\varepsilon})+\log\,\varepsilon\right)+C_2\leq c_{\varepsilon}^{\frac{1}{n}}\tau_j\log(|z|+\varepsilon)+C_2.
			$$
			Letting $\varepsilon\rightarrow0$, we get that $\varphi\leq c^{\frac{1}{n}}\tau_j\log|z|+O(1)$ almost everywhere and hence everywhere since both sides lie in $\mbox{PSH}(U)$.
			
		\end{proof}
		
		\subsection{The extended Demailly-P\u{a}un conjecture}
		The main purpose of this section is to give a partial answer to an  extended version of the  Demailly-P\u{a}un's conjecture \cite[Conjectures  1.2]{LWZ24}, extending the result in \cite[Theorem 1.11]{LWZ24}.
		
		\begin{theorem}\label{thm:DP conjecture}
			Let $(X,\omega)$ be a compact Hermitian manifold. Let $\{\beta\}\in BC^{1,1}(X)$ be a Bott-Chern class satisfying $\underline{\mbox{Vol}}(\beta)>0$ and $\overline{\mbox{Vol}}_{n-1}(\beta)<+\infty$, assume also that $\beta$ has a bounded potential $\rho\in \mbox{PSH}(X,\beta)\cap L^{\infty}(X)$. Then the class $\{\beta\}$ is big, i.e., it contains a Hermitian current.
		\end{theorem}
		
		We first make some preparations.

		\begin{definition}
			A function $u$ is said to be $(\beta+dd^c\rho)$-psh if $u+\rho$ is $\beta$-psh. Therefore, we have a one-to-one correspondence between $\mbox{PSH}(X,\beta)$ and $\mbox{PSH}(X,\beta+dd^c\rho)$.
		\end{definition}
		\begin{theorem}\label{thm:DP lemma 1}
			Let $(X,\omega,\{\beta\})$ be as in \cref{thm:DP conjecture}. Let $0\leq f\in L^p(X,\omega^n)$ be such that $\int_Xf\omega^n>0$ for some $p>1$. Let $\lambda>0$ be a constant. Then, there exists a $\varphi_\lambda\in \mbox{PSH}(X,\omega+\beta+dd^c\rho)$ with $\rho+\varphi_{\lambda}\in C(X)$, satisfying:
			$$
			(\omega+\beta+dd^c\rho+dd^c\varphi_{\lambda})^n=e^{\lambda\varphi_{\lambda}}f\omega^n
			$$
			in the weak sense of currents.
		\end{theorem}
		\begin{proof}
			The proof can be copied along the lines of the proof in \cite[Theorem 6.4]{LWZ24}, the closeness of the reference form $\beta$ does not seem to be required.
		\end{proof}
		
		The following criterion is due to Lamari \cite{Lam99}.
		\begin{lemma}\label{lem:criterion} (\cite{Lam99}). Let $\alpha$ be a smooth real $(1,1)$-form. Then, there exists a distribution $\psi$ on $X$ such that $\alpha+dd^c\psi\geq0$ if and only if:
			$$
			\int_X\alpha\wedge\gamma^{n-1}\geq0
			$$    
			for any Gauduchon metric $\gamma$ on $X$.
		\end{lemma}
		
		\begin{lemma}\label{lem: DP conjecture 2}
			Let $g$ be a Gauduchon mertic on $X$ and set $G:=g^{n-1}$. For $0<\varepsilon<1$, let $v_{\varepsilon}\in \mbox{PSH}(X,\beta+dd^c\rho)\cap L^{\infty}(X)$ be the unique solution to 
			$$
			(\beta+dd^c\rho+dd^cv_{\varepsilon})^n=e^{\varepsilon v_{\varepsilon}}\omega\wedge G.
			$$
			The existence and uniqueness follows from \cref{MA equation 1}. Then we have
			$$
			\left(\int_X(\beta+dd^c\rho+dd^cv_{\varepsilon})\wedge G\right)\left(\int_X(\beta+dd^c\rho+dd^cv_{\varepsilon})^{n-1}\wedge\omega\right)\geq\frac{1}{n}\left(\int_Xe^{\frac{\varepsilon v_\varepsilon}{2}}\omega\wedge G\right)^2.
			$$
		\end{lemma}
		
		\begin{proof}
			The proof can be copied without any change from \cite[Lemma 6.7]{LWZ24}, the closeness of $\beta$ is not required.
		\end{proof}
		
		\begin{proof}[Proof of \cref{thm:DP conjecture}]
			We argue by contradiction. Suppose $\{\beta\}$ does not adimit any Hermitian current. Then by \cref{lem:criterion}, there is a sequence of positive numbers $\delta_j$ decreases to $0$ and Gauduchon metrics $g_j$ such that
			\begin{equation}\label{eq 1}
				\int_X(\beta+dd^c\rho-\delta_j\omega)\wedge g_j^{n-1}\leq0.
			\end{equation}
			Set $G_j:=g_j^{n-1}$. Then the above inequality reads
			\begin{equation}\label{eq 2}
				\int_X(\beta+dd^c\rho)\wedge G_j\leq\delta_j\int_X\omega\wedge G_j.
			\end{equation}

			\textbf{Claim}. For each $j$ there is a function $v_j\in PSH(X,\beta+dd^c\rho)\cap L^{\infty}(X)$ solves 
			\begin{equation}\label{eq 0}
				(\beta+dd^c\rho+dd^cv_j)^n=c_j\omega\wedge G_j,\quad \sup_Xv_j=0
			\end{equation}
			where $c_j>0$ is a constant, and satisfies
			$\left(\int_X\beta_j\wedge G_j\right)\left(\int_X\beta_j^{n-1}\wedge\omega\right)\geq\frac{c_j}{n}\left(\int_X\omega\wedge G_j\right)^2$
			for $\beta_j:=\beta+dd^c\rho+dd^cv_j$.
			
			\vspace{\baselineskip}
			
			For clarity, we defer the proof of the claim to later.
			
			Integrating both sides of (\ref{eq 0}) we get
			\begin{equation}\label{eq 3}
				c_j\geq\frac{\underline{\mbox{Vol}}(\beta)}{\int_X\omega\wedge G_j}>0.
			\end{equation}
			Recall that $\beta_j:=\beta+dd^c\rho+dd^cv_j$. By (\ref{eq 1}) we have
			\begin{equation}\label{eq 4}
				\int_X\beta_j\wedge G_j=\int_X(\beta+dd^c\rho)\wedge G_j\leq\delta_j\int_X\omega\wedge G_j.
			\end{equation}
			
			combining (\ref{eq 2}) and the above claim we have
			$$
			\frac{c_j}{n}\left(\int_X\omega\wedge G_j\right)^2\leq\left(\int_X\beta_j\wedge G_j\right)\left(\int_X\beta_j^{n-1}\wedge\omega\right)\leq\left(\delta_j\int_X\omega\wedge G_j\right)\left(\int_X\beta_j^{n-1}\wedge\omega\right)
			$$
			and hence $\frac{c_j}{n}\int_X\omega\wedge G_j\leq\delta_j\int_X\beta_j^{n-1}\wedge\omega$. Combining this inequality with (\ref{eq 3}) we get
			$$
			0<\frac{1}{n}\underline{\mbox{Vol}}(\beta)\leq\frac{c_j}{n}\int_X\omega\wedge G_j\leq\delta_j\int_X\beta_j^{n-1}\wedge\omega\leq\delta_j\overline{\mbox{Vol}}_{n-1}(\beta).
			$$
			The RHS $\rightarrow0$ because $\overline{\mbox{Vol}}_{n-1}(\beta)<+\infty$ which leads to a contradiction. The proof is therefore concluded.
		\end{proof}
		\begin{proof}[Proof of the claim]
			By the Bedford-Taylor convergence theorem \cite{BT82}, it suffices to prove that $\rho+v_{\varepsilon}-M_{\varepsilon}\rightarrow v_j$ uniformly on $X$ for some constants $M_{\varepsilon}$ and $e^{\varepsilon v_{\varepsilon}}\rightarrow c_j$ for some positive constant $c_j$ as $\varepsilon\rightarrow0$. Where we take $G=G_j$  in \cref{lem: DP conjecture 2} and $v_\varepsilon$ are the corresponding solutions. Considering the Monge-Amp\`ere equations:
			\begin{equation}\label{eq 6}
				(\beta+dd^c(\rho+v_\varepsilon-M_{\varepsilon}))^n=e^{\varepsilon(\rho+v_{\varepsilon}-M_{\varepsilon})+\varepsilon M_{\varepsilon}}(e^{-\varepsilon\rho}\omega\wedge G_j),
			\end{equation}
			where $M_{\varepsilon}:=\sup_X(\rho+v_{\varepsilon})$. 
			
			We first observe that $\varepsilon M_{\varepsilon}$ are uniformly bounded. Indeed, since $\rho$ is bounded, we may assume $\rho\geq0$,  integrating both sides of (\ref{eq 6}) we have
			$$
			e^{\varepsilon M_{\varepsilon}}\int_X\omega\wedge G_j\geq e^{\varepsilon M_{\varepsilon}}\int_Xe^{-\varepsilon\rho}\omega\wedge G_j\geq\underline{\mbox{Vol}}(\beta)>0
			$$
			and hence
			$$
			\varepsilon M_{\varepsilon}\geq \log\frac{\underline{\mbox{Vol}}(\beta)}{\int_X\omega\wedge G_j},
			$$
			which gives the uniform lower bound of $\varepsilon M_{\varepsilon}$. For the upper bound, set $\psi_\epsilon:=\rho+v_\epsilon-M_\epsilon\in\mbox{PSH}(X,\beta)$ with $\sup_X\psi_\epsilon=0$ and $h_\epsilon:=e^{-\epsilon\rho}\omega\wedge G_j/\omega^n$, which is a uniformly bounded (away from zero) family of positive functions on $X$. Then \eqref{eq 6} can be rewritten as
			$$
			(\beta+dd^c\psi_\epsilon)^n=e^{\epsilon M_\epsilon}\cdot e^{\epsilon\psi_\epsilon}\cdot h_\epsilon\omega^n.
			$$
			The mixed type inequality (\cref{mixed type 2}) gives that
			$$
			(\beta+dd^c\psi_\epsilon)\wedge\omega^{n-1}\geq e^{\frac{\epsilon M_\epsilon}{n}}\cdot (h_\epsilon e^{\epsilon \psi_\epsilon})^{\frac{1}{n}}\omega^n.
			$$
			Integrating both sides we have
			\begin{align*}
				e^{\frac{\epsilon M_\epsilon}{n}}\int_Xh_\epsilon^{\frac{1}{n}}\cdot e^{\frac{\epsilon\psi_\epsilon}{n}}\omega^n\leq\int_X\beta\wedge\omega^{n-1}+\int_X\psi_\epsilon dd^c(\omega^{n-1})\leq C,
			\end{align*}
			where the last inequality is because $\sup_X\psi_\epsilon=0$ and hence $\psi_\epsilon$ lies in a compact subset of $\mbox{PSH}(X,\beta)$. To obtain an upper bound for $\epsilon M_\epsilon$, it suffices to obtain a lower bound of $\int_Xh_\epsilon^{\frac{1}{n}}\cdot e^{\frac{\epsilon\psi_\epsilon}{n}}\omega^n$. To see this, we can assume without loss of generality that $\omega^n(X)=1$ and using Jensen's inequality to get
			$$
			\int_Xh_\epsilon^{\frac{1}{n}}\cdot e^{\frac{\epsilon\psi_\epsilon}{n}}\omega^n\geq \inf_Xh_\epsilon^{\frac{1}{n}}\cdot e^{\int_X\frac{\epsilon}{n}\psi_\epsilon}\omega^n\geq c_0>0,
			$$
			the last inequality is again due to the fact that $\psi_\epsilon$ lies in a compact subset of $\mbox{PSH}(X,\beta)$ and hence we derive the upper bound of $\epsilon M_\epsilon$.
			
			Note that the LHS of (\ref{eq 6}) is uniformly bounded in $L^p(X,\omega^n)$ since $\varepsilon M_{\varepsilon}$ is uniformly bounded. By \cref{thm: a prior estimate} we have $\psi_{\varepsilon}:=\rho+v_\varepsilon-M_{\varepsilon}$ is uniformly bounded and $\sup_X\psi_{\varepsilon}=0$. We may thus assume without loss of generality that $\varepsilon M_{\varepsilon}\rightarrow M$ and $\psi_{\varepsilon}$ is a Cauchy sequence in $L^1(X,\omega^n)$. Since $(\beta
			+dd^c\psi_{\varepsilon})^n\leq e^{\varepsilon M_{\varepsilon}-\varepsilon\rho}\omega\wedge G_j\leq\lambda\omega^n$ for some uniform constant $\lambda>0$, we can use \cref{cor:stability 2} to deduce that
			$$
			||\psi_\varepsilon-\psi_{\delta}||_{\infty}\leq||\psi_\varepsilon-\psi_{\delta}||^\gamma_{L^1(X,\lambda\omega^n)}.
			$$
			Consequently, $\psi_{\varepsilon}$ is Cauchy in $L^{\infty}(X)$ and we may assume that $\psi_{\varepsilon}$ converges uniformly to $\psi$ and $\sup_X\psi=0$. Therefore, by the Bedford-Taylor convergence theorem \cite{BT82} we have
			$$
			(\beta+dd^c(\rho+v_\varepsilon-M_{\varepsilon}))^n\rightarrow(\beta+dd^c\psi)^n=e^M\omega\wedge G_j.
			$$
			We can thus take $c_j:=e^M$ and $v_j:=\psi-\rho-\sup_X(\psi-\rho)$ and the claim follows.

		\end{proof}

	\end{document}